\documentclass{conm-p-l}
\usepackage{amsmath,amssymb}
\usepackage[all]{xy}

\theoremstyle{plain}

\newtheorem{theorem}{Theorem}[section]

\newtheorem{problem}[theorem]{Problem}

\theoremstyle{definition}

\theoremstyle{remark}

\def\G{{\mathcal G}}
\def\R{{\mathbb R}}
\def\Z{{\mathbb Z}}
\def\Q{{\mathbb Q}}
\def\C{{\mathbb C}}
\def\L{{\mathbb L}}
\def\M{{\mathcal M}}
\def\A{{\mathcal A}}
\def\LL{\mathcal{L}}
\def\H{{\mathbb H}}
\def\aut{\mathit{aut}}
\def\Hol{\mathrm{Hol}}
\def\ad{\mathrm{ad}}
\def\Ad{\mathrm{Ad}}
\def\Lo{\mathit{L}}
\def\map{\mathit{map}} 
\def\Der{\mathrm{Der}}
\def\Hom{\mathrm{Hom}}
\def\Rel{\mathrm{Rel}}
\def\cat0{\mathrm{cat}_0}

\def\dim{\mathrm{dim}}
\def\ker{\mathrm{ker}}

\def\im{\mathrm{im}}

 \makeatletter
\def\@strippedMR{}
\def\@scanforMR#1#2#3\endscan{%
  \ifx#1M\ifx#2R\def\@strippedMR{#3}%
  \else\def\@strippedMR{#1#2#3}%
  \fi\fi}
\renewcommand\MR[1]{\relax
  \ifhmode\unskip\spacefactor3000 \space\fi
  \@scanforMR#1\endscan
  MR\MRhref{\@strippedMR}{\@strippedMR}}
\makeatother 

\begin{document}

\title[The Homotopy Theory of Function Spaces: A Survey]
{The Homotopy Theory of Function Spaces: A Survey}

\author{Samuel Bruce Smith}

\address{Department of Mathematics,
   Saint Joseph's University,
   Philadelphia, PA 19131}

\email{smith@sju.edu}

\date{\today}

\keywords{Function space, monoid of self-equivalence, free loop space, space of holomorphic maps, gauge group, string topology,  configuration space,  section space, classifying space, Gottlieb group, localization, rational homotopy theory}

\subjclass[2010]{ 55P15, 55P35,  55P48, 55P50, 55P60, 55P62, 55Q52, 55R35, 46M20}

\begin{abstract} We  survey  research  on the homotopy theory of the space  $\map(X, Y)$ consisting of all continuous functions between two topological spaces.   
We  summarize  progress on various classification problems for   the homotopy types represented by  the path-components  of $\map(X, Y)$. We  also discuss work  on the homotopy theory of the monoid of self-equivalences $\aut(X)$ and of the free loop space $\Lo X$.  We     consider these topics in both ordinary  homotopy theory as well as after localization. In the latter case, we discuss  algebraic models     for the localization of  function spaces and  their applications. 
\end{abstract}

\maketitle
\section{Introduction.}
In this paper,  we survey research in homotopy theory on   function spaces treated as topological spaces  of interest in their own right.  
 We begin, in this section,  with some general remarks on the topology of function spaces.  We then give a  brief historical sketch of work on the homotopy theory of function spaces. This sketch   serves to introduce the basic themes   around which the  body of the paper is   organized.

By  work of {\sc Brown} \cite[1964]{MR0165497} and {\sc Steenrod}  \cite[1967]{MR0210075},   the homotopy theory of    function spaces  may be studied    in the ``convenient  category'' of compactly generated Hausdorff spaces.  Retopologizing   is required, however.   Given spaces $X$ and $Y$ in this category, let $Y^X$ 
 denote the space of all continuous functions   with the compact-open topology.  Define  $$\map(X, Y) = \emph{k}\,(Y^X)$$ to be the associated compactly generated space.  Then   $\map(X, Y)$   satisfies the    desired exponential laws and is a homotopy
 invariant of $X$ and $Y$. The space  $\map(X, Y)$ is generally disconnected with path-components corresponding to the set  of free homotopy classes of maps.   We write $\map(X, Y;f )$ for the path-component containing a given map $f \colon X \to Y.$ Important special cases include:      $\map(X, Y;0),$   the space of null-homotopic maps;   $\map(X,X;1),$  the identity component; $\aut(X),$ 
 the space of all homotopy self-equivalences of $X$; and $\Lo X = \map(S^1, X)$ the free loop space. 
 
Concrete results on the path-components of  $\map(X,Y)$   often require much more restrictive
hypotheses on $X$ and $Y.$  By  {\sc Milnor} \cite[1959]{MR0100267}, when $X$ is a compact, metric space and $Y$ is a CW complex, the components $\map(X, Y;f)$ are   of CW homotopy   type. 
A natural case  to consider then is when $X$ is a finite CW complex and $Y$ is any CW complex.  
By {\sc Kahn} \cite[1984]{MR733413},  $\map(X, Y)$ is also of CW type when $X$ is any CW complex and $Y$ has finitely many homotopy groups.  

The space $\map(X, Y)$ has two close relatives. If $X$ and $Y$ are pointed spaces, we have $\map_*(X, Y)$   the space of basepoint preserving functions, with components $\map_*(X, Y; f)$  for $f$ a based map.   Given a fibration $p \colon E \to X,$   we have $\Gamma(p)$ the space of sections with components $\Gamma(p; s)$ for $s$ a fixed section. Of course, $\map(X, Y) \simeq \Gamma(p)$ when $p$ fibre-homotopy trivial with fibre $Y$.    
 Many theorems about $\map(X, Y)$ generalize to $\Gamma(p)$  and many have related  versions for $\map_*(X, Y)$.   For the sake of  brevity,  when possible we  state   theorems for the free function space and omit  extensions and restrictions.    Theorems stated for the based function space   are then results that do not apply to $\map(X, Y)$.

  \subsection{A Brief History.}
Function spaces are at the foundations of homotopy theory  and appear in the literature dating back, at least,  to Hurewicz's definition of the homotopy groups in the 1930s. Work focusing   explicitly on the  homotopy theory of a function space first       appears in the 1940s.  
 {\sc Whitehead}  \cite[1946]{MR0016672}  introduced the problem of 
classifying the homotopy types represented by the path-components of a function space, focusing on the case $\map(S^2, S^2)$.    {\sc Hu}  \cite[1946]{MR0019920} showed $$ \pi_1(\map(S^2, S^2; \iota_m)) \cong \Z/2|m| ,$$
where $\iota_m$ is the map of degree $m$ thus distinguishing components of different absolute degree.  

A decade later,        papers of {\sc Thom} \cite[1957]{MR0089408} and {\sc Federer} \cite[1956]{MR0079265}  appeared giving dual methods for  computing homotopy groups of components of $\map(X, Y)$.  Thom        used a Postnikov decomposition of $Y$ to indicate a method of calculation.   Federer constructed a spectral sequence converging to these homotopy groups using a cellular decomposition of $X$.   Both authors obtained the following basic identity:
$$ \label{Thom} 
\pi_q(\map(X, K(\pi, n); 0)) \cong H^{n-q}(X; \pi)$$
for $X$ a CW complex and $\pi$ an abelian group.  

In the 1960s, the  monoid $\aut(X)$ of all homotopy self-equivalences of $X$  emerged as  a central object for the theory of fibrations.  {\sc Stasheff} \cite[1963]{MR0154286} constructed a universal fibration  for  CW  fibrations   with  fibre of the homotopy type  of a fixed finite CW complex $X,$ building on work of {\sc Dold-Lashof} \cite[1957]{MR0101521}.        His result implied  the universal $X$-fibration is obtained, up to homotopy,  by applying the  Dold-Lashof  classifying space functor   to the evaluation fibration $\omega \colon \map(X, X; 1) \to X$. 
In this same period, {\sc Gottlieb} \cite[1965]{MR0189027}  introduced and studied  the  evaluation subgroups or {\em Gottlieb  groups}:   $$G_n(X) = \im \{ \omega_\sharp \colon \pi_n(\map(X, X; 1)) \to \pi_n(X) \} \subseteq \pi_n(X)$$  initiating a rich literature 
on the evaluation map.  Among many other properties, he showed
the Gottlieb groups correspond to the image of     the linking homomorphism in the long exact sequence of homotopy groups of   the universal $X$-fibration.  Thus  the vanishing of a Gottlieb group $G_n(X)$   is equivalent to  the vanishing of  the linking homomorphism in  degree $n$ for every CW fibration with fibre $X$.

In the   1970s,    {\sc Hansen} \cite[1974]{MR0362396} began a systematic study of the homotopy classification problem for the path components of $\map(X, Y)$. He completed the classification for  $\map(S^n, S^n)$   building   on the  methods    of Whitehead, mentioned above.    He and other authors   obtained      complete   results  in many special  cases   involving   spheres, suspensions, projective spaces and  certain manifolds.   

The space of holomorphic maps $\Hol(M, N)$ between two complex manifolds offers a  deep refinement of the homotopy  classification  problem for continuous maps  with important  interdisciplinary connections. {\sc Segal} \cite[1979]{MR533892}  initiated the study of the space  
  $\Hol(M, N)$ in homotopy theory  proving  the inclusion  $$\Hol^*_k(S^2, \C P^m) \hookrightarrow \map_*(S^2, \C P^m; \iota_k)$$ 
 induces a homology equivalence through a range of degrees. Here $\Hol^*_k(S^2, \C P^m)$ denotes the space of based holomorphic maps of degree $k$.  In fundamental work in complex geometry,  {\sc Gromov} \cite[1989]{MR1001851} identified the class of  elliptic     manifolds and proved they satisfy  the ``Oka Principle''.  As a consequence, he identified a large class  of manifolds    for which the inclusion 
 $\Hol(M, N) \hookrightarrow \map(M, N)$ is a weak equivalence.   {\sc Cohen-Cohen-Mann-Milgram} \cite[1991]{MR1097023} described the full   stable homotopy type of  $\Hol^*_k(S^2, \C P^m),$ their description  given in terms of  configuration spaces.  A related problem of stabilization for moduli spaces of connections is the subject of the famous ``Atiyah-Jones conjecture'' in mathematical physics  {\sc Atiyah-Jones} \cite[1978]{MR503187}.  

The gauge groups provide a  connection between the homotopy   theory of function spaces  and the theory of principal bundles.    
   Let $P \colon E \to X$ be a principal $G$-bundle
for $G$ a connected topological group classified by a map $h \colon X \to BG$.
The {\em gauge group} $\G(P)$ of $P$ is   defined to be the group of $G$-equivariant
homeomorphisms $f \colon E \to E$ over $X$.     {\sc Atiyah-Bott} \cite[1983]{MR702806} used the      gauge group in their celebrated study of Yang-Mills  equations and principal bundles over a Riemann surface. They made use of Thom's theory and a multiplicative equivalence originally due to     {\sc Gottlieb}  \cite[1972]{MR0309111}    $$ \G(P) \simeq \Omega \map(X, BG; h)$$
to study the homotopy theory of  $B\G(P)$.  Gottlieb's identity, in turn,  links 
the classification of gauge groups up to $H$-homotopy type, for fixed $G$ and $X,$   to the homotopy classification problem for   $\map(X, BG)$.  {\sc Crabb-Sutherland} \cite[2000]{MR1781154}
proved that the gauge groups $\G(P)$ represent only finitely many homotopy types  for $G$ a  compact Lie group  and $X$ a finite complex,. In  contrast,  the  path-components of $\map(X, BG)$ may represent infinitely many distinct homotopy types in this case by {\sc Masbaum} \cite[1991]{MR1101938}.

The advent of localization techniques introduced new depth  to the study of function spaces while opening up a wide range of fundamental problems and applications.
In his seminal paper on rational homotopy theory,  {\sc Sullivan} \cite[1977]{MR0646078}  sketched a 
construction for an algebraic model for components of $\map(X, Y)$ for $X$ and $Y$ simply connected CW complexes with $X$ finite, as an extension of Thom's ideas. Sullivan's construction
was completed by {\sc Haefliger} \cite[1982]{MR667163}. Sullivan also identified the rational Samelson Lie algebra of $\aut(X)$ for $X$ a finite, simply connected CW complex via an isomorphism:
$$ \pi_*(\aut(X)) \otimes \Q,   [ \,\,  , \, \, ] \, \cong  \, H_*(\Der(\mathcal{M}_X)), [ \, \,  , \, \,   ]$$
Here the latter space is the homology of the Lie algebra of degree lowering derivations 
of the Sullivan minimal model of $X$ with the commutator bracket. 
 
 One of the early  applications of Sullivan's rational homotopy theory was  the proof by {\sc Vigu{\'e}-Poirrier-Sullivan} \cite[1976]{MR0455028}   of the unboundedness of the  Betti numbers of   the free loop space
 $\Lo X = \map(S^1, X)$ for certain   simply connected CW complexes $X$.
 Combined with a famous result of  
 {\sc   Gromoll-Meyer} \cite[1969]{MR0264551} in geometry,  this calculation solved the  ``closed geodesic problem'' for many manifolds.
     The calculation was    deduced from a  Sullivan model constructed   for $\Lo X.$    
 
 The $p$-local homotopy theory of a function space featured in  a landmark result in algebraic topology,   the proof of the Sullivan conjecture.   {\sc Miller} \cite [1984]{MR750716} proved  $$\pi_n(\map_*(B\pi, X; 0)) = 0 \, \, \, \hbox{\, for all \,} n \geq 0$$ 
where  $\pi$ is any finite group and   $X$ any finite CW complex. Among many applications, this result was used by {\sc McGibbon-Neisendorfer} \cite[1984]{MR749108} to affirm Serre's 
conjecture:   $\pi_m(X)$ contains a subgroup of order $p$ for  infinitely many  $m$.

 {\sc Lannes} \cite[1987]{MR932261} constructed the $T$-functor which is    left adjoint to the tensor product in the category of unstable modules over the Steenrod algebras. His construction provided a model for the mod $p$ cohomology of the space $\map(BV, X)$  where $V$ is a $p$-group.   Lannes' construction was adapted to the rational homotopy setting by {\sc Bousfield-Peterson-Smith} \cite[1989]{MR989883} and, later,  {\sc Brown-Szczarba} \cite[1997]{MR1407482}  to give another model
 for the rational homotopy type of $\map(X, Y;f)$. {\sc Fresse} \cite[preprint]{Fressepreprint}  recently constructed 
 a version of Lannes' functor in a category of operadic algebras giving a model for the integral 
 homotopy type of certain function spaces.

 The free loop space recently re-emerged as a central object for study in homotopy theory with the appearance of  work of {\sc Chas-Sullivan} \cite[preprint]{preprintCS}. They   constructed a   product on the regraded homology  
 $$\mathbb{H}_*(\Lo M^m) = H_{*+m}(\Lo M^m)$$ for a simply connected, closed, oriented $m$-manifold $M^m$ using intersection theory. They also defined a bracket  on the equivariant  homology of $\Lo M^m$ and a degree $+1$ operator  giving $\mathbb{H}_*(\Lo M^m)$ the structure of Batalin-Vilkovisky algebra. These    structures have   incarnations in diverse other settings.  Their study,  known as {\em string topology}, is now an active subfield in the intersection of  homotopy theory and geometry. 
 
\subsection{Organization.}
In Section \ref{sec:classification},   we discuss work  on  the ordinary  and stable homotopy theory, as opposed to the local homotopy theory  of   function spaces.   We focus on the areas introduced above, namely: (i) the  general path component $\map(X, Y; f);$     (ii) the   monoid $\aut(X);$  and (iii)  the free loop space $\Lo X.$ We also discuss work on the stable homotopy theory of these spaces and on   spectral sequence calculations of their invariants. 
    In Section \ref{sec:local}, we discuss   the localization of function spaces.
      We describe the algebraic models of Sullivan, and of later authors, for the general component,  the monoid of self-equivalences and the free loop space in rational  homotopy theory,  and survey their applications. We also discuss  the $p$-local homotopy theory of function spaces including the work of Miller, Lannes and others on  the space of maps out of  a classifying space,  and algebraic models for function spaces in tame homotopy theory. The paper includes a rather extensive bibliography  gathering together  both  papers directly focused on function spaces and papers giving significant applications and extensions.   


%

\section{Ordinary and Stable Homotopy Theory of Function Spaces.}
\label{sec:classification}
We divide our discussion in this section according to the    cases (i), (ii) and (iii)  above.  We then discuss some general results in  stable homotopy theory and spectral sequence constructions for function spaces.

\subsection{General  Components.}

As mentioned in the introduction, the following open problem lies at the historical roots of the study of function spaces as objects in their own right. 

 \begin{problem} \label{prob:class} Given spaces $X$ and $Y$ classify the path-components $\map(X, Y; f)$
up to homotopy type for homotopy classes  $f \colon X \to Y.$
\end{problem}
 
  We consider a variety of cases here beginning with the most classical, mentioning progress  on Problem \ref{prob:class},  when appropriate.

\subsubsection{Maps from Spheres and Suspensions. }

The components of $\map(S^p, Y)$ correspond to the homotopy classes in  $\pi_p(Y)$. The coproduct on $S^p$ gives rise to an equivalence
 $\map_*(S^p, Y; \alpha ) \simeq \map_*(S^p, Y; 0),$ for any class $\alpha \colon S^p \to Y$.
  By adjointness,    
$\pi_n(\map_*(S^p, Y; 0))  \cong \pi_{n+p}(Y).$   These observations were made by {\sc Whitehead} \cite[1946]{MR0016672} who gave the first algebraic method for computation. Whitehead identified   the linking homomorphism in the long exact homotopy sequence for the evaluation fibration $\map_*(S^p, Y; \alpha) \to \map(S^p, Y;  \alpha) \to Y$
obtaining:
$$ \xymatrix{ \cdots \ar[r] & \pi_{n+1}(Y) \ar[r]^{\! \! \! \! \! \! \!  \! \! \!  \! \! \!  \partial}  \ar[rd]_{W(\alpha)} & \pi_n(\map_*(S^p, Y; \alpha)) \ar[d]^\cong \ar[r] & \pi_n(\map(S^p, Y; \alpha)) \ar[r] & \cdots \\
& & \pi_{n+p}(Y) & & }
$$
where $W(\alpha)(\beta) = -[\alpha, \beta]_w$ denotes the Whitehead product map. 
Using this sequence, he proved   $\map(S^2, S^2; \iota) \not\simeq \map(S^2, S^2; 0)$  by comparing homotopy groups.   {\sc Hu} \cite[1946]{MR0019920} and {\sc Koh} \cite[1960]{MR0119201} computed $\pi_{2m-1}(\map(S^{2m}, S^{2m}; \alpha))$  for small values of $m$.  In these cases, the order  of $\pi_{2m-1}(\map(S^{2m}, S^{2m}; \alpha))$ depends on the absolute value of the degree of $\alpha$ and so distinguishes components with different absolute order. Since clearly $\map(S^{2m}, S^{2m}; \alpha) \simeq \map(S^{2m}, S^{2m}; -\alpha)$ the classification in these cases was complete with these calculations.

{\sc Hansen} \cite[1974]{MR0362396} obtained the complete  classification for self-maps of  $S^n$. For even spheres,   he proved 
$$\map(S^{2m}, S^{2m}; \alpha) \simeq \map(S^{2m}, S^{2m}; \beta) \iff  [\alpha, \iota]_w = \pm [\beta, \iota]_w. $$   
Here $\iota \in \pi_{2m}(S^{2m})$ is the fundamental class. For odd spheres,   the components of $\map(S^{2m-1}, S^{2m-1})$ are all homotopy equivalent for $m = 1, 2, 4$ due to the existence of a multiplication on $S^{2m-1}$ in these cases. For   $m \neq 1, 2, 4$,  Hansen showed   $\map(S^{2m-1}, S^{2m-1}; \iota) \not\simeq \map(S^{2m-1}, S^{2m-1}; 0)$ 
and 
$$\map(S^{2m-1}, S^{2m-1}; \alpha) \simeq  \left\{ \begin{array}{l}  \map(S^{2m-1}, S^{2m-1}; \iota)  \hbox{ \, if  deg}(\alpha) = \hbox{odd} \\
  \map(S^{2m-1}, S^{2m-1}; 0)  \hbox{ \, if  deg}(\alpha) = \hbox{even}. \end{array} \right. $$
 
 Problem \ref{prob:class}  remains open for $\map(S^{m}, S^{n})$ for $m > n$. 
 {\sc Yoon} \cite[1995]{MR1369453}  observed a connection between the Gottlieb group $G_m(Y)$ and the    homotopy classification problem  for  $\map(S^m, Y) $ showing $\map(S^m, Y; \alpha) \simeq \map(S^m, Y; 0)$ if and only if $\alpha \in G_m(Y).$     {\sc Lupton-Smith} \cite[2008]{MR2431641}
extended this to a surjection of sets
$$ \xymatrix{ \pi_m(Y)/ G_m(Y) \ar@{->>}[r] &
\begin{displaystyle}
\frac{\{ \text{components   $ 
\map(S^m,Y;f)$} \}}{\text{homotopy equivalence}}
\end{displaystyle}.}
$$
Thus the complexity of the classification problem for $\map(S^{m}, S^{n})$
is roughly that of computing Gottlieb groups $G_m(S^n)$. 
   Extensive, low-dimensional calculations of this group were recently made by {\sc Golasi{\' n}ski-Mukai}  \cite[2009]{MR2427733}.   {\sc Lee-Mimura-Woo} \cite[2004]{MR2100869} calculated the Gottlieb groups for certain homogeneous spaces.

When $X = \Sigma A$ is a suspension, the fibres $\map_*(\Sigma A, Y; f)$ of the various  evaluation fibrations $\omega_f \colon \map(\Sigma A, Y; f) \to Y$ are all homotopy equivalent to the space $\map_*(\Sigma A, Y; 0)$ with homotopy groups  
$$ \pi_q(\map_*(\Sigma A, Y; 0)) = [\Sigma^{q+1}A, Y].$$  
{\sc Lang} \cite[1973]{MR0341484} extended Whitehead's exact sequence to  this case replacing the Whitehead product in $\pi_*(Y)$ by the generalized Whitehead product in $[\Sigma^* A, Y].$ 
It is natural to consider, as 
Whitehead did,  a   stronger version of Problem \ref{prob:class}, namely, the  classification of  the evaluation fibrations $ \omega_f \colon \map(X, Y; f) \to Y$
up to  fibre homotopy type for homotopy classes $f\colon  X \to Y$.
   {\sc Hansen} \cite[1974]{MR0368000} defined    
$\omega_f \colon \map(\Sigma A, \Sigma B; f) \to Y$  to be {\em strongly fibre homotopy equivalent}
to  $\omega_g \colon \map(\Sigma A, \Sigma B; g) \to \Sigma B$ if the fibre homotopy equivalence is homotopic to the identity after (fixed) identification of the fibres with    $\map_*(\Sigma A, \Sigma B; 0).$
He proved:
$$ \omega_f \hbox{\, is strongly fibre homotopic to \,} \omega_g \iff [f, 1_{\Sigma B}] = [g, 1_{\Sigma B}]$$
where $[ \, , \, ]$ here denotes the generalized Whitehead product in $[\Sigma A, \Sigma B]$. 
{\sc McClendon} \cite[1981]{MR640979} showed that the evaluation fibrations $\omega_f \colon \map(\Sigma A, Y;f) \to Y$ behave as principal fibrations and, in particular,  are   classified by maps $s \colon Y \to \map(A, Y)$
  determined by generalized Whitehead products.  
    
 \subsubsection{Maps into Eilenberg-Mac\! Lane Spaces.}
 
The weak  homotopy type of  the space  $\map(X, K(\pi, n); f)$ may be described for any $f \colon X \to K(\pi, n)$ for $\pi$ abelian.  The ideas are due to {\sc Thom} \cite[1957]{MR0089408} with a refinement by {\sc Haefliger} \cite[1982]{MR667163}.  
 First, observe that these components are all homotopy equivalent since $K(\pi, n)$ has the homotopy type of a topological group. A homotopy class $\alpha \in \pi_p(\map(X, K(\pi, n; 0)))$ corresponds, by adjointness, to a map $A \colon S^p \times X \to K(\pi, n)$.   On cohomology,   $$A^*(x_n) = 1 \otimes a_n + u_p \otimes a_{n-p}$$ 
 where $a_n, a_{n-p} \in H^*(X; \pi)$ with subscripts indicating degree while $u_p \in H^p(S^p; \pi)$ and $x_{n} \in H^n(K(\pi, n); \pi)$ are the fundamental classes.
Since $A$ restricts to the constant map on   $S^p \times *$ we see $a_n = 0$. The assignment $ \alpha \mapsto a_{n-p}$   gives the identification $$\pi_p(\map(X, K(\pi, n); f)) \cong H^{n-p}(X; \pi),$$ 
mentioned in the introduction and leads to directly to a  
weak equivalence 
$$ \map(X, K(\pi, n); f)) \simeq_w \prod_{p \geq 1} K(H^{n-p}(X; \pi), p).$$
Thom  also indicated how   the homotopy groups $\pi_p(\map(X, Y; f))$ for $Y$ a finite Postnikov piece are determined, up to extensions,  by the $k$-invariants of $Y$ and the groups $H^{n-p}(X, \pi_n(Y)).$ This approach was encoded in Haefliger's construction of a Sullivan model for $\map(X, Y;f)$, as discussed below.

{\sc Gottlieb} \cite[1969]{MR0239595}  extended Thom's result to the case $n =1$ and $\pi$ any group. Here
$$\map(X, K(\pi, 1); f) \simeq_w K(C(f_\sharp), 1)$$
where $C(f_\sharp)$ denotes the centralizer of the image of $f_\sharp \colon \pi_1(X) \to \pi.$ 
{\sc M{\o}ller} \cite[1987]{MR910659} showed that when $Y$ is a twisted
Eilenberg-Mac Lane space, then $\map(X, Y; f)$ is one also with homotopy groups determined  by the  cohomology groups of $X$ with twisted coefficients in the homotopy
groups of $Y$.  Note that the weak equivalences above are homotopy equivalences by Whitehead's Theorem, when $\map(X, K(\pi, n))$ is of CW type, e.g., when $X$ is compact or a CW complex. 
 In general, the study of the homotopy type of $\map(X, Y ;f)$ when $Y$ has at least two nonvanishing  homotopy 
groups is a difficult, open problem. 

\subsubsection{Maps between Manifolds.}
  The homotopy theory of $\map(M^m, N^n)$ for $M^m$ and $N^n$ closed manifolds
 is a topic of wide-ranging interest.  In this case, important variations  have been considered. Below we consider one such variation with direct ties to Problem \ref{prob:class}, namely      spaces of  holomorphic maps.  
 We begin with the space $\map(M^m, N^n)$.
 
 If $T_g$ is an orientable surface, then $T_g \simeq K(\pi_1(T_g), 1)$ and the classification problem for $\map(X, T_g)$ reduces to the computation of centralizers of homomorphisms into $\pi_1(T_g)$. For $g \geq 2$ this group is  highly nonabelian and the only possible nontrivial centralizers are isomorphic to $\Z$ by {\sc Hansen} \cite[1983]{MR703761}.
{\sc Hansen}  \cite[1974]{MR0385919} earlier considered the space $\map(T_g, S^2)$.  As a generalization of Whitehead's exact sequence, he showed
an exact sequence
$$ 0 \to \Z/2|m| \to \pi_1(\map(T_g, S^2; \iota_m)) \to \Z^{2g} \to 0 $$
 which gives the classification, in terms of degree,  in this case.  
 The   fundamental group $\pi_1(\map(T_g, S^2; \iota_m))$ was later completely determined by 
 {\sc Larmore-Thomas} \cite[1980]{MR612697}.  
 
 {\sc Hansen} \cite[1981]{MR607120} extended his classification result for the space of self-maps of spheres to the case
 $\map(M^m, S^m)$ where $M^m$ is  closed, oriented, connected manifold with vanishing first Betti number. Note that, by Hopf's Theorem, $[M^m, S^m] = \Z$ with maps classified by degree. When $m$ is even and $\geq 4, $ Hansen showed the homotopy types of 
 $\map(M^m, S^m; \alpha)$ are classified by the absolute values of the degrees of the $\alpha.$ When $m$ is odd and 
 $m \neq 1,4,7$ there are two homotopy types corresponding to the distinct types $\map(M^m, S^m; 0)$
 and $\map(M^m, S^m; \iota)$ where $\iota$ is of degree $1.$  Again in this case, components are classified by the parity of degree of the class $\alpha$.
 
  {\sc Sutherland} \cite[1983]{MR698208}  extended Hansen's work eliminating the restriction on the first Betti number and dealing with the  case $M^m$ nonorientable. Note that in the latter case there are only
  two distinct classes $\alpha \colon M^m \to S^m$ and so the problem reduces to distinguishing these components for $m \neq 1,4, 7$. Sutherland observed that the components of $\map(M^m, S^m)$ all have the same homotopy type if   there
  is a map $I \colon M^m \to \map(S^m, S^m; \iota)$ making the diagram
  $$ \xymatrix{ & \map(S^m, S^m; 1) \ar[d]^\omega \\ M^m \ar[ur]^{I} \ar[r]^\iota & S^m}$$
   commute, where  $\iota$  is of degree $1$. Taking $M^m = \R P^m$,
   we have a lift  $I$ based on the lift  $I' \colon \R P^m \to SO(n+1) \subseteq \map(S^m, S^m; 1)$ of $\iota.$ Sutherland showed the components of $\map(M^m, S^m)$
   are all of the same homotopy type if there exists a map $f \colon M^m \to \R P^m$
   of odd degree giving examples with $m \neq 1, 4, 7$ for which all the components
   are homotopy equivalent.

   {\sc Sasao} \cite[1974]{MR0346783} studied the homotopy type of components of $\map(\C P^m, \C P^n; i)$
   for $m \leq n$ and $i \colon \C P^m \to \C P^n$ the inclusion. He constructed  a map
   $$ \alpha_{m, n}  \colon U(n+1)/ \Delta(m+1) \times U(n-m) \to \map(\C P^m, \C P^n; i)$$
  where $\Delta(m+1) \subset U(m+1)$ denotes scalar multiplies of the identity.  He proved $\alpha_{m,n}$ induces an isomorphism on rational homotopy groups and on ordinary homotopy groups through degree $4n-4m+1$.  {\sc Yamaguchi} \cite[1983]{MR717319} extended Sasao's analysis to the case of quaternionic projective spaces. 
   {\sc M{\o}ller} \cite[1984]{MR744651} gave the complete classification for the components $\map(\C P^m, \C P^n)$
   showing the homotopy types are classified by the absolute value of the degree of 
   a representative class.  The result is a direct consequence of his calculation
   $$H_{2n-2m+1}(\map(\C P^m,\C P^n; \iota_k)) \cong \Z/d \hbox{\, where \,} d =  \binom{n+1}{m} | k|^m. $$
 {\sc Yamaguchi}   \cite[2006]{MR2028664} studied maps between real projective spaces.  He defined the analogue of Sasao's   map, here of the form 
  $$\beta_{m,n}  \colon  O(n+1)/\Delta(m+1) \times O(n-m) \to \map(\R P^m,\R P^n; i)$$
 and proved $\beta_{m,n}$ is an equivalence on ordinary and rational homotopy groups through certain ranges of degrees. 
 
 When $G$ is a topological group (or group-like space) the path-components of $\map(X, G)$ are all of the same homotopy type.  Problem \ref{prob:class}  thus reduces, in this case, to the study of the homotopy theory of the null-component $\map(X, G; 0)$.  Given Lie groups $G$ and $H$, the calculation of homotopy invariants of $\map(G, H)$  is a difficult open problem.       Recently,  {\sc Maruyama-{\=O}shima} \cite[2008]{MR2440413} computed the homotopy groups of $\map_*(G, G)$ for $G = SU(3), Sp(2)$ in degrees $\leq 8.$

\subsubsection{Spaces of Holomorphic Maps.}
{\sc Segal} \cite[1979]{MR533892} proved  a  basic result on  the homotopy theory of the space $\Hol(M, N)$. His work launched a vital subfield of research on the   ``stability'' of  the inclusion $\Hol(M, N) \hookrightarrow \map(M, N).$ 
 Segal proved  $$\Hol^*_k(T_g, \C P^n) \hookrightarrow \map_*(T_g, ; \iota_k)$$  induces a homology  isomorphism  through dimension $(k-2g)(2n-1)$ where $T_g$ is a Riemann surface of genus $g$. Specializing to the case of the sphere, he proved $$ \Hol^*_k(S^2, \C P^n) \hookrightarrow \map_*(S^2, \C P^n; \iota_k) $$  
induces a homotopy equivalence up to degree $2n -1.$

Segal's work was extended by many authors. {\sc Guest} \cite[1984]{MR729753} 
proved the corresponding stability result on homology for  $\Hol_k^*(S^2, F) \hookrightarrow \map(S^2, F; \iota_k)$ for certain complex flag manifolds $F$.  His proof involved developing  the analogue of a  Morse-theoretic result for the case of  the  energy functional on the space $C^\infty(S^2, F)$ of smooth maps. 
{\sc Kirwan} \cite[1986]{MR884188}  extended Segal's result to the case the target is the complex Grassmannian manifold $G(n, n+m)$ of $n$-planes in $n+m$-space   proving $\Hol_k^*(S^2, G(n, n+m)) \hookrightarrow \map(S^2,G(n, n+m); \iota_k)$ induces a homology isomorphism in degrees depending on $k, n$ and $m.$  
{\sc Mann-Milgram} \cite[1991]{MR1094457}
considered this case as well, constructing a spectral sequence to analyze the homology of 
$\Hol_k^*(S^2, G(n, n+m)).$ 
{\sc Gravesen} \cite[1989]{MR989398} studied holomorphic maps into space $\Omega G$ for $G$ a complex, compact Lie group.

{\sc Cohen-Cohen-Mann-Milgram} \cite[1991]{MR1097023}   and {\sc Cohen-Shimamota} \cite[1991]{MR1113390} described  the stable homotopy type of $\Hol_k^*(S^2, \C P^n)$. They proved $$\Hol^*_k(S^2, \C P^n) \simeq C_k(\R^2, S^{2n-1})$$
where $C_k(\R^2, S^{2n-1})$ is the configuration space of distinct points in $\R^2$
with labels in $S^{2n-1}$ of length at most $k$.   
 Cohen-Cohen-Mann-Milgram also computed the homology of $\Hol_k^*(S^2, \C P^n)$ with $\Z_p$-coefficients in terms of Dyer-Lashof operations. 
{\sc Mann-Milgram} \cite{MR1231702} used the stable homotopy decomposition
above to  prove the homology stability of the inclusion $\Hol_k^*(S^2, F) \hookrightarrow \map_*(S^2, F)$ for $F$ an $SL(n, \C)$-flag-manifold. 
{\sc Boyer-Mann-Hurtubise-Milgram} \cite[1994]{MR1294670} and  {\sc Hurtubuise} 
\cite[1996]{MR1661615}  proved    homology stabilization theorems for the  space $\Hol_k^*(S^2, G/P)$  for certain complex homogeneous spaces $G/P$.  {\sc Boyer-Hurtubise-Milgram} gave a configuration space description of $\Hol_k(T_g, M)$ for certain complex manifolds admitting nice  Lie group actions extending the approach  of Gravesen.  

Segal's stabilization problem  has deep interdisciplinary connections. 
 {\sc Gromov} \cite[1989]{MR1001851} obtained general stability results as a consequence of  his work on the  {\em Oka Principle} in complex geometry. A complex manifold $M$ satisfies the Oka principle if every continuous map $f \colon S \to M$ is homotopic to a holomorphic map where $S$ is a Stein manifold.  Gromov  identified  the class of  ``elliptic'' manifolds and proved elliptic manifolds satisfied the Oka principle.  Consequently, he obtained the inclusion $$\Hol(S, M) \hookrightarrow  \map(S, M)$$ is a weak homotopy equivalence for $S$   Stein  and $M$  elliptic. The class of elliptic manifolds includes complex Lie groups and their homogeneous spaces.     
 
 The problem of stabilization also has a famous   incarnation in Yang-Mills theory and mathematical physics.   {\sc Atiyah-Jones} \cite[1978]{MR503187} constructed a
 map $$\theta_k \colon \mathcal{M}_k \to  \map_*(S^3, SU(2); \iota_k )$$
 where $\mathcal{M}_k$ is a moduli space of connections  on  a principal $SU(2)$-bundle $P_k$ over $S^4$ corresponding to a map $S^4 \to BSU(2)$ of degree $k$. 
 They proved $\theta_k$ induces a homology surjection through a range of degrees and conjectured $\theta_k$ induces an equivalence in both homotopy and homology through a range depending on $k$. Work on the Atiyah-Jones  conjecture includes 
 {\sc Taubes} \cite[1989]{MR978084},  {\sc Gravesen} \cite[1989]{MR989398} and  {\sc Boyer-Hurtubise--Mann-Milgram} \cite[1993]{MR1217348}.
  
 Many authors have studied related  spaces of  maps. 
  {\sc Vassiliev} \cite[1992]{MR1089669}  proved a stable equivalence
$$ \Hol^*_k(S^2, \C P^n) \simeq \mathrm{SP}_{n-1}^k(\C)$$
where the latter is the space of monic complex polynomials of degree $k$ with all roots 
of multiplicity $< n.$ 
{\sc Guest-Kozlowski-Yamaguchi} \cite[1994]{MR1306670} extended Segal's result in a 
different direction, proving the inclusion 
$$\Hol_k^*(S^2, X_n) \hookrightarrow \map_*(S^2, X_n)$$ is a homotopy equivalence up to degree $k$ where $X_n \subseteq \C P^{n-1}$ is the subspace of points with at most
one coordinate zero.  
The cohomology of the space of basepoint-free holomorphic maps  $\Hol_1(S^2, S^2)$ was studied by {\sc Havlicek} \cite[1995]{MR1318152} while the homotopy groups of  $\Hol_k(S^2, S^2)$ were studied by {\sc Guest-Kozlowski-Murayama-Yamaguchi} \cite[1995]{MR1365252}. 
{\sc Kallel-Milgram} \cite[1997]{MR1601616} gave a complete calculation of the homology of 
$\Hol^*_k(T_g, \C P^m)$ for $T_g$ an elliptic  Riemann surface. 
 The space of real rational functions was recently studied by {\sc Kamiyama} \cite[2007]{MR2316447}.

 {\sc Kallel-Salvatore} \cite[2006]{MR2284046}  applied techniques from string topology to the study of spaces of maps between manifolds.  Set $$\mathbb{H}_*(\map(M^m, N^n)) = H_{*+n}( \map(M^m, N^n))$$ and similarly for $\Hol(S^2, N^n).$ When $M^m, N^n$ are closed, compact and orientable, they proved $\mathbb{H}_*(\map(S^m, N))$ has   a ring structure corresponding to an intersection product and  
 $\mathbb{H}_*(\map(M^m, N^n))$ is a module over this ring. 
They used this structure to compute  $\mathbb{H}_*(\map(S^2, \C P^n; \iota_k))$ and $\mathbb{H}_*(\Hol_k(S^2, \C P^n))$ with   $\Z_p$-coefficients.    They also studied $\mathbb{H}_*(\map(T_g, \C P^n; \iota_k); \Z_p))$ for $T_g$ a compact Riemann surface proving, among other results, that these groups are   isomorphic for all $k$ when $ p$  divides $n.$

 \subsubsection{Maps into a Classifying Space and Gauge Groups.}
 
Let $X$ be a space and $G$ a connected topological group. Suppose $P \colon E \to X$ is a principal $G$-bundle.   The gauge group $\mathcal{G}(P)$ is the topological group of bundle automorphisms  of $P$.      The gauge group featured in important work of  
 {\sc Atiyah-Bott} \cite[1983]{MR702806} in mathematical physics.  They considered the action of $\G(P)$ on the moduli space $\mathcal{A}$ of Yang Mills connections on a principal $U(n)$-bundle $P \colon E \to M$ for $M$ a  Riemann manifold.    Among other results, they proved  $H^*(B\G(P))$ is torsion free and computed its    Poincar{\' e} series.  
 Their calculation depends on the  identity: 
     $$ \G(P) \simeq_H \Omega \map(X, BG; h),$$
     where $h \colon X \to BG$ is the classifying map of $P$, a result
 originally due to {\sc Gottlieb}  \cite[1972]{MR0309111}.  Here $X$ is a finite CW complex. 
Thus $B\G(P) \simeq \map(X, BG; h).$  By Bott periodicity, the loops and double loops   on $BU(n)$ are torsion free. Atiyah-Bott  used this fact and Thom's theory
to make their  calculations.

The classification of gauge groups for fixed $X$ and $G$ up to $H$-equivalence
or, alternately, up to ordinary equivalence is the subject of active research. 
Gottlieb's identity implies  the homotopy classification problem for $\map(X, BG)$  refines   the gauge group classification problem.  
{\sc Masbaum} \cite[1991]{MR1101938} studied the homology of the components of the space $\map(X, BSU(2))$ for 
$X$ a $4$-dimensional CW complex obtained by attaching a single $4$-cell to a bouquet
of $2$-spheres. This case includes oriented, simply connected $4$-dimensional manifolds.  Using a cofibre sequence for $X$, he obtained, in particular, that the components of $\map(S^4, BSU(2))$ represent 
infinitely many homotopy types.  Using a related analysis, {\sc Sutherland} \cite[1992]{MR698208} considered  the  classification of components of $\map(T_g, BU(n))$ for $T_g$ an orientable  surface of genus $g$. He obtained the calculatiuon
$$ \pi_{2n-1}(\map(T_g, BU(n); \iota_k)) \cong \Z^g \oplus \Z/d \hbox{ \, where \,} d = (n-1)!(k, n)$$  thus distinguishing the components corresponding to maps $k$ and $l$  with $(k, n) \neq (l, n).$ {\sc Tsukuda} \cite[2001]{MR1834088} classified the homotopy types represented by the components  of $\map(S^4, BSU(2))$ showing   $\map(S^4, BSU(2); \iota_k) \simeq  \map(S^4, BSU(2); \iota_l)$ if and only if $ k = \pm l$. {\sc Kono-Tsukuda} \cite[2000]{MR1381542} generalized this result from $X= S^4$ to $X$ a simply connected $4$-dimensional manifold.
  
As regards the homotopy type of the gauge group, {\sc Kono} \cite[1991]{MR1103296} proved $$\G(P_k) \simeq \G(P_l) \iff (12, k) = (12,  l)$$ for $k,l \in \pi_4(BSU(2)) \cong \Z.$ Thus the infinitely many distinct homotopy types represented by the path-components of  $\map(S^4, BSU(2))$ loop to only $6$ distinct homotopy types.   Kono's proof  included    the calculation $$\pi_2(\G(P_k)) = \Z/ (12, k)$$      using, essentially,  Whitehead's exact sequence mentioned above.      {\sc Kono-Tsukuda}  \cite[1996]{MR1381542} extended this result from $X = S^4$ to $X$ a closed, simply connected manifold using a cofibre sequence for $X$ to make the corresponding calculation.  {\sc Hamanaka-Kono} \cite[2006]{MR2310471} obtained  a corresponding classification for $SU(3)$-bundles over $S^4$.    They proved $$\G(P_k) \simeq \G(P_l) \iff (120, k) = (120, l)$$
where  $P_k$ and $P_l$ are principal $SU(3)$-bundles over $S^6$ with third Chern class equal to $2k$ and $2l$, respectively. 
 
 {\sc Crabb-Sutherland} \cite[2000]{MR1781154} obtained a global result on the classification of gauge groups. They proved  that, for any fixed finite CW complex $X$ and  compact Lie group $G$, there are only finitely many $H$-homotopy types represented by the gauge groups $\G(P)$ for $P$ a principal $G$-bundle over $X.$ 
  Their proof is based on an alternate description of the gauge group: 
 $$ \G(P) \cong \Gamma(\Ad(P))$$ where $\Ad(P) = E \times_G G^{\ad} \to X$ is the adjoint bundle associated to $P \colon E \to X$ and  the space of sections has   multiplication induced by $G$. A key step in their finiteness result is the proof that the fibrewise rationalization of $\Ad(P)$ is equivariantly trivial. They also classified the $H$-homotopy types of gauge groups of  $SU(2)$-principal bundles over $S^4$ complementing Masbaum and Kono's work. Here
 $$\G(P_k) \simeq_H \G(P_l) \iff (180, k) = (180, l).$$
Summarizing,   the   infinitely many distinct homotopy types represented by the components
of $\map(S^4, BSU(2))$ loop to $6$ distinct homotopy types and 
   $18$ distinct  $H$-homotopy types.  
 
 \subsection{Spaces of Self-Equivalences.} The space of equivalences $\aut(X)$  of a space $X$ with some additional structure admits many important refinements.   When $M$ is a   Riemannian manifold, we have the chain of subspaces
 $$\mathrm{Isom}(M) \hookrightarrow   \mathrm{Diff}(M) \hookrightarrow  \mathrm{Homeo}(M) 
 \hookrightarrow  \aut(M)$$ 
 given by spaces of isometries, diffeomorphisms and homeomorphisms, respectively. 
 Each of these spaces is the subject of active research in homotopy theory and geometric topology.    
 {\sc Smale} \cite[1959]{MR0112149}   proved the inclusion $$\mathrm{Isom}(S^2) \hookrightarrow \mathrm{Diff}(S^2)$$ is a homotopy equivalence. Since $\mathrm{Isom}(S^2) \simeq O(3)$ this determines the homotopy type of $\mathrm{Diff}(S^2).$ The (Generalized) Smale Conjecture asserts that $\mathrm{Isom}(M) \simeq \mathrm{Diff}(M)$  for $M$ a $3$-manifold of constant, positive curvature.  The Smale Conjecture was affirmed for $M = S^3$ by {\sc Hatcher} \cite[1983]{MR701256}.   {\sc Gabai} \cite[2003]{MR1895350} proved the corresponding result for $M$ a closed, hyperbolic $3$-dimensional manifold.

The first result on the $H$-homotopy type of $\aut(X)$ is due, essentially, to {\sc Thom} \cite[1957]{MR0089408}. By   his results mentioned above, we have
$$ \aut(K(\pi, n)) \simeq_H \aut(\pi) \rtimes K(\pi, n)$$
for $\pi$ abelian.  When $n =1$, Gottlieb's extension of Thom's results leads to an identification
$$ \aut(K(\pi, 1)) \simeq_H \mathrm{Out}(\pi) \rtimes K(\pi, 1)$$
where $\mathrm{Out}(\pi)$ denotes the group of outer automorphisms. 
Note  that these results include a description of    $\pi_0(\aut(X))$,  the group   of free homotopy self-equivalences of $X$.  This group is, in general,  quite complicated even for simple $X$.    See   {\sc Arkowitz} \cite[1990]{MR1070585} and {\sc Rutter} \cite[1997]{MR1474967}  for surveys of the extensive literature on this group.

 Since the path-components of $\aut(X)$ are all of the same homotopy type, we focus on     the  component
of the identity which we denote $\aut(X)_\circ$.   Thus
$$\aut(X)_\circ = \map(X, X ;1)$$ 
is the identity component in the space of self-maps.  

{\sc Hansen} \cite[1990]{MR1070574} identified the homotopy type of $\aut(S^2)_\circ$ by comparing the evaluation fibration for this space with the fibre sequence $SO(2) \to SO(3) \to S^2.$
He proved
$$ \aut(S^2)_\circ  \simeq_H SO(3) \times \widetilde{\aut_*(S^2)_\circ}$$
where $\widetilde{Z}$ denotes the universal cover. Combined with Smale's result, this shows $\mathrm{Diff}(S^2) \hookrightarrow \aut(S^2)$ is not a homotopy equivalence. 
{\sc Yamanoshita} \cite[1993]{MR1219881} obtained a related result  proving
$$ \aut(\R P^2)_\circ   \simeq_H SO(3) \times \left(\aut_*(\R P^2)_\circ) /O(2)\right).$$
This result implies $\mathrm{Diff}(\R P^2) \simeq O(3)$ is not  homotopy equivalent
to $\aut(\R P^2).$ {\sc Yamanoshita} \cite[1985]{MR792987}  also obtained a general result 
$$ \aut(X \times Y) \simeq \aut(X) \times \aut(Y) \times \map_*(Y, \aut(X)) \times \map_*(X, \aut(Y))$$  provided the dimension of $X$ is less than the connectivity of $Y$. In particular:
$$ \aut(S^1 \times Y) \simeq O(2) \times \aut(Y) \times \Omega \aut_\circ(Y)$$
 for simply connected $Y.$
 
{\sc McCullough} \cite[1981]{MR597873} computed $\pi_q(\aut(M)_\circ)$ for $1 \leq q \leq n-3$ for $M$ a connected sum of closed, aspherical manifolds of dimension $\leq 3$ proving the  groups $\pi_{n-2}(\aut(M)_\circ )$  are not finitely generated. He used this result to give examples of closed $3$-manifolds $M$ such that the fundamental group of  $\mathrm{Homeo}(M)$ is not finitely generated.
{\sc Didierjean} \cite[1990]{MR1070573} and  \cite[1992]{MR986023}  used a Postnikov decomposition of $X$ to construct a spectral sequence  converging to the  homotopy groups of $\aut(X)$. She determined low degree homotopy groups of $\aut(X)$ for  $X = SU(3)$ and $X = Sp(2)$ up to extensions.

Given a fibration $p \colon E \to B,$ we may consider the monoid  $\aut(p)$ of fibre-homotopy equivalences $f \colon E \to E$ covering the identity of $B.$  {\sc Booth-Heath-Morgan-Piccinini} \cite[1984]{MR743373}   extended Gottlieb's result for the gauge group  to prove an $H$-equivalence
$$ \aut(p) \simeq_H \Omega \map(B, B\aut(F); h))$$ where
$F$ is the fibre of $p$ and $h \colon B \to B\aut(F)$ is the classifying map. A simplicial version of  this result was earlier obtained by {\sc Dror-Dwyer-Kan} \cite[1980]{MR581012}. 
 {\sc Didierjean} \cite[1987]{MR810666} extended Thom's result to the fibrewise setting, proving
$$\pi_q(\aut(p)) \cong H^{n-q}(B; \pi)$$ for $p \colon E \to B$ a principal fibration with fibre $K(\pi, n)$ and   made calculations of the homotopy groups of $\aut(p)$
when $F$ has two nonzero homotopy groups. She constructed a  spectral sequence converging to the homotopy groups of $\aut(p)$,  expanding on work of   {\sc Legrand} \cite[1983]{MR732851}. 

\subsection{The Free Loop Space.} The space of maps
$$\Lo X = \map(S^1, X)$$  is the subject of intensive  research in
 diverse branches of mathematics.  Given a compact Lie group $G$, the space of smooth loops on $G$  is 
 an infinite-dimensional Lie group, called    the  ``loop group'' of $G$.    The representation theory of loop groups    has  deep connections to mathematical physics (cf. {\sc Pressley-Segal} \cite[1986]{MR900587}).     {\sc Gromoll-Meyer} \cite[1969]{MR0264551}  linked the closed geodesic problem for  a  compact Riemannian manifold $M$ to the homotopy theory of the free-loop space $\Lo M.$
They proved $M$ admits  infinitely many closed, prime geodesics   if     the  Betti
numbers of $\Lo M$ grow without bound.   {\sc Vigu{\'e}-Poirrier-Sullivan} \cite[1976]{MR0455028}  proved that  the Betti numbers of $\Lo M$  are unbounded when the rational cohomology of $M$ requires at least two generators.  Their calculation was facilitated by  a Sullivan model for $\Lo M$,  described in Section \ref{sec:local}, below.  

More recently,  {\sc Chas-Sullivan} \cite[preprint]{preprintCS} and \cite[2004]{MR2077595} unearthed a wealth of structure on the (regraded) homology of $\Lo M$ of a   closed, oriented, smooth $m$-manifold $M.$  Setting $ \H_*(\Lo M) = H_{*+m}(\Lo M)$  they defined a graded-commutative and associative product
   $\bullet$ on $\H_*(\Lo M)$     and a  related Lie bracket  on  the equivariant  homology. The pair 
give   $\H_*(\Lo M)$ the structure of Gerstenhaber algebra. A
degree $+1$ operator $\Delta$   give $\H_*(\Lo M)$ 
the structure of Batalin-Vilkovisky (BV) algebra.  These structures were obtained by geometric methods using intersection theory  and  transversality arguments.

A homotopy theoretic construction of the string topology  structures was given by {\sc Cohen-Jones} \cite[2002]{MR1942249} using Thom spectra.   They also proved, for $M$ simply connected,  an isomorphism of graded algebras
$$\H_*(\Lo M; \mathbb{F}) \cong HH^*(S^*(M), S^*(M); \mathbb{F})$$ 
where the latter  is the Hochschild cohomology of the algebra of singular cochains
of $M.$ Here $\mathbb{F}$ is a field. 
The Cohen-Jones construction was extended to more general ring spectra by 
{\sc Gruher-Salvatore} \cite[2008]{MR2392316}. {\sc Cohen-Klein-Sullivan}
\cite[2008]{MR2399136},   {\sc Crabb} \cite[2008]{MR2391632}, and  Gruher-Salvatore  independently proved the homotopy invariance of the loop product and bracket, a significant advance since the original constructions depended on the smooth structure of $M.$    {\sc Chataur} \cite[2005]{MR2180465},   {\sc Hu} \cite[2006]{MR2251161} and  
 {\sc Kallel-Salvatore} \cite[2006]{MR2284046}  considered   generalizations of string operations  from $\Lo M = \map(S^1, M)$ to $\map(S^n, M)$.  The work of these various authors include   constructions   of the string topology operations  in the frameworks of ring spectra (Cohen-Jones  and Gruher-Salvatore]), bordism theory (Chataur) and fibrewise homotopy theory (Crabb). 

  As regards  the ordinary   homotopy theory of the free loop space, {\sc Hansen} \cite[1974]{MR0341477} gave an example of an aspherical space $X$ with $\pi_1(\Lo X)$ not finitely generated.  Note that when $X$ is an $H$-space $\Lo X \simeq X \times \Omega X$.  {\sc Aguad{\' e}} \cite[1987]{MR915024} made a general study of spaces $X$, called {\em $T$-spaces},  for which the
  evaluation fibration $$\Omega X \to \Lo X \to X$$ is fibre-homotopically trivial. He obtained a refinement of the notion of $H$-space via a sequence of classes $T=T_1 \subset T_2 \subset  \ldots \subset T_\infty = H$-spaces with separating examples.  {\sc Woo-Yoon} \cite[1995]{MR1346627} proved that when $X$ is a $T$-space the components of $\map(\Sigma A, X)$ are all homotopy equivalent. {\sc Fadell-Husseini} \cite[1989]{MR984789}  proved $\Lo  X$ has infinite L.S. category
  for $X$ a simply connected CW complex with finitely generated, nontrivial rational cohomology.

   {\sc Smith}  \cite[1981]{MR630771} and \cite[1984]{MR748960} constructed an Eilenberg-Moore spectral sequence for the cohomology of a free loop space.
   Starting with the pull-back square
$$ \xymatrix{ \Lo X \ar[r]^{\omega} \ar[d]_{\omega} & X \ar[d]^{  \Delta} \\
X \ar[r]^{\! \! \! \! \! \! \! \! \Delta} & X \times X}$$ he obtained a  spectral sequence
$$E^{*,*}_2 =  \mathrm{Tor}_{H^*(X \times X; \mathbb{F})}(H^*(X; \mathbb{F}), H^*(X; \mathbb{F})) \implies 
H^*(\Lo X; \mathbb{F}).$$
He proved  collapsing results for this spectral sequence  and obtained  calculations of 
$H^*(\Lo M; \mathbb{F})$ for $\mathbb{F}$ of characteristic $2$ and $0$. 
 {\sc Kuribayashi} \cite[1991]{MR1096437} used the Eilenberg-Moore spectral sequence   to prove  the fibre  is totally noncohomologous to zero in the fibration $\Omega M \to \Lo M \to M$   for $M$ a Grassmann or Stiefel manifold and $\mod p$ cohomology   for certain primes $p$. 
 
{\sc McCleary-Ziller} \cite[1987]{MR900038}  proved  the Betti numbers of  $\Lo M$ are unbounded 
for $M$ a compact, simply connected homogeneous space  not equivalent to a rank one symmetric space using spectral sequence methods and extending earlier work of {\sc Ziller} \cite[1977]{MR0649625} who used Morse theory. {\sc Roos} \cite[1988]{MR981826} studied the Poincar{\' e}-Betti series for $\Lo X$ for $X$ a wedge of  spheres using local algebra. He proved the series for $X = S^2 \vee S^2$ is not
rational. {\sc Halperin-Vigu{\' e}-Poirrier} \cite[1991]{MR1084712} proved the  $\mathbb{F}$-Betti numbers are unbounded for  a field $\mathbb{F}$ of positive characteristic $k$  provided $H^*(X; \mathbb{F})$ requires at least $2$ generators and under certain restrictions on $k$ and $X.$  
{\sc McCleary-McLaughlin} \cite[1992]{MR1079897} studied the free loop space in the context of Morava $K$-theory while {\sc Ottosen} \cite[2003]{MR2009581} considered the Borel cohomology of the free loop space. 
{\sc Lambrechts} \cite[2001]{MR1834084}  proved the Betti numbers of the free loop space are unbounded for certain connected sums.

 {\sc Burghelea-Fiedorowicz} \cite[1984]{MR750675} and {\sc Goodwillie} \cite[1985]{MR793184}  proved an isomorphism
 of graded spaces $$H^*(\Lo  X) \cong HH^*(S_*(\Omega X), S_*(\Omega X))$$
 where the latter space is the Hochschild cohomology. {\sc Menichi} \cite[2001]{MR1854644}, {\sc Dupont-Hess} \cite[2002]{MR1918186} and {\sc Ndombol-Thomas} \cite[2002]{MR1871242}  independently proved this is an isomorphism of algebras. Menichi also made
 calculations of the graded algebra  $H^*(\Lo  X; \Z_p)$ for $X$ a suspension and $X = \C P^n$ while Ndombol-Thomas
  used Hochschild cohomology to make calculations of  $H^*(\Lo  X; \Z_p)$ for $X = S^m, \C P^m$ and $\Sigma \C P^m$.
  {\sc Kuribayashi-Yamaguchi} \cite[1997]{MR1472851} made complete calculations of $H^*(\Lo X; \Z_p)$ when $X$ is simply connected with mod $p$ cohomology an exterior algebra on few generators using  Hochschild cohomology to obtain  information on the $E_2$-term in the Eilenberg-Moore spectral sequence.
    Recently, {\sc Seeliger} \cite[2008]{MR2370366} used the Serre spectral sequence applied to the evaluation fibration to make
  calculations of $H^*(\Lo  \C P^m).$

  Since the appearance of the paper of Chas-Sullivan, the structure of $\H_*(\Lo M)$
  has seen an explosion of research with many partial descriptions of the loop product, the  Gerstenhaber algebra and the BV algebra structure  in special cases.  We mention a sampling of these results here. {\sc Cohen-Jones-Yan} \cite[2004]{MR2039760} constructed a spectral sequence 
of algebras $$E^2_{*,*} = H^*(M; H_*(\Omega M)) \Longrightarrow \H_*(\Lo M)$$ 
and used this to calculate the loop product for $M = S^2, \C P^n$.  {\sc Tamanoi} \cite[2006]{MR2211159}  computed the BV algebra structure of $\H_*(\Lo M)$ for $M = SU(n)$ and $M$ a complex Stiefel manifold.   {\sc Gruher-Salvatore} \cite[2008]{MR2392316} extended the string operations to the case $M = BG$ for $G$ a compact Lie group. {\sc Menichi} \cite[2009]{MR2466078}  proved the Cohen-Jones isomorphism $$\H_*(\Lo S^2; \mathbb{F}_2) \cong HH^*(S^*(S^2;  \mathbb{F}_2), S^*(S^2;  \mathbb{F}_2))$$ mentioned above is an isomorphism of Gerstenhaber algebras but, surprisingly, not an isomorphism of BV algebras.

 \subsection{Spectral Sequences and Stable Decompositions.}
 We here discuss some general results on homotopy invariants and the stable homotopy type of function spaces not covered by the preceding discussion.
 
{\sc Federer} \cite[1956]{MR0079265} constructed a spectral sequence converging to the homotopy groups of $\map(X, Y; f)$ for $X$ any finite CW complex and $Y$ a simple CW complex.  He defined an exact couple  from the long exact homotopy sequences  of the restriction
 maps $\rho_n \colon \map(X_{n}, Y; f_n) \to \map(X_{n-1}, Y; f_{n-1})$. He  identified the homotopy groups of the fibre of $\rho_n$ with cellular cochain groups of $X$ with coefficients in $\pi_*(Y)$ and obtained a spectral sequence   $$ E_2^{p,q} = H^q(X; \pi_{p+q}(Y)) \implies \pi_p(\map(X, Y; f)).  $$  
{\sc Dyer} \cite[1966]{MR0203728} applied the Federer spectral sequence to calculate  low degree homotopy and homology group of components of $\map(X, Y)$ when dim$X$ is less than the connectivity of $Y.$ {\sc Schultz} \cite[1973]{MR0407866} and {\sc M{\o}ller} \cite[1990]{MR1038731}  constructed   equivariant  versions of this spectral sequence.  
  
   {\sc Borsuk} \cite[1952]{MR0056285} proved   that if $X$ is a finite CW complex of dimension $k$ with nonzero $k$th Betti number then $\map_*(X, S^n)$ has
 nonzero $(n-k)$th Betti number. {\sc Moore} \cite[1956]{MR0083123} extended this analysis and asked for a spectral sequence relating the cohomology of $X$ to the homology of $\map_*(X, S^n)$. 
 {\sc Spanier} \cite[1959]{MR0105106} showed the functor $$F_n(X) = \varinjlim_k \map(X, \Omega^k \mathrm{SP}^{n+k} S^{n+k}),$$
 converts cohomology groups to homotopy groups. Here  $SP^{n+k}$ is the symmetric product functor. 
  {\sc Anderson} \cite[1972]{MR0310889} constructed an Eilenberg-Moore spectral sequence for the
  cohomology  of   $\map_*(X, Y)$ using  the cobar construction in the category of cosimplicial spaces.   {\sc Legrand} \cite[1986]{MR925276}   constructed   a spectral sequence in the spirit of Moore's problem for $\map_*(X, Y)$ using a Postnikov decomposition of $Y$. 
{\sc Patras-Thomas} \cite[2003]{MR1956614} proved the Anderson spectral sequence converges when the dimension of $X$ is no bigger than the connectivity of $Y.$  {\sc Chataur-Thomas} \cite[2004]{MR2100872} 
gave a related $E_\infty$-model   for function spaces. Applied to the free loop space, their model gives an operadic version of Hochschild cohomology. 

The stable homotopy type of the based function space $\map_*(X, Y)$ has been described in many cases.  When $X= S^n$, this is   just the loop space $\Omega^n Y$. Stable decompositions of this space are   central structural results in homotopy theory. {\sc Snaith} \cite[1974]{MR0339155} proved 
a stable decomposition for $\Omega^n Y$   in terms of configuration spaces of $j$-tuples of distinct points of $\R^n$.  His proof was obtained by analyzing the stable homotopy type of approximations due to {\sc May} \cite[1972]{MR0420610} for  iterated loop spaces.
{\sc B{\" o}digheimer}  \cite[1987]{MR922926} proved a generalization of Snaith's  splitting result. He showed
$$\Omega^\infty \Sigma^\infty \map_*(M, \Sigma^n Y) \simeq \Omega^\infty \Sigma^\infty \left(C(M, \partial M; n) \wedge_{\Sigma_n} Y^n\right) $$ where $M$ is a compact manifold with boundary and $ C(M, \partial M; n)$ is the configuration space. 
{\sc Arone} \cite[1999]{MR1638238} described the Goodwillie tower of  $\Omega^\infty \Sigma^\infty \map_*(X, \Sigma^n Y)$ for more general $X$ recovering the previous splittings. He obtained  an Eilenberg-Moore spectral sequence from this description. {\sc Ahearn-Kuhn} \cite[2002]{MR1917068} and {\sc Kuhn} \cite[2006]{MR2196061} studied this spectral sequence proving it a spectral sequence of graded algebras and studying functorial properties.

{\sc Campbell-Cohen-Peterson-Selick} \cite[1987]{MR921473}  studied   $\map_* (P_m(2^r),S^n)$ where $P_m(2^r) =S^{m-1}\cup e^m$ is a Moore space with   attaching map   of degree $2^r$. They  gave a partial description of   the mod $2$ Steenrod operations and  proved that  $\map_*(P_3(2), S^n)$     is not  decomposable as a product  except, perhaps,  for   finitely many $n$. 
 {\sc Westerland} \cite[2006]{MR2248387} obtained a stable splitting for components of the space  $\map_*(T_g,S^2)$ where $T_g$ is a surface of genus $g > 0.$
 {\sc Cohen} \cite[1987]{MR922928} and  {\sc B{\"o}digheimer} \cite[1987]{MR922926} independently obtained a stable splitting of the free loop space $\Lo  \Sigma X$  in terms of configuration spaces.   Other results on the stable homotopy of the free loop space include   splitting results for $\Lo  \R P^n$ by {\sc  Bauer-Crabb-Spreafico}  \cite[2001]{MR1817001} and {\sc Yamaguchi}  \cite[2005]{MR2134058}.

\section{Localization of   Function Spaces.}
\label{sec:local}
In this section, we survey work on function spaces after localization.  Recall a {\em nilpotent space} $X$ is a connected CW complex   such that 
$\pi_1(X)$ is a nilpotent group and the standard action of $\pi_1(X)$ on the higher homotopy groups of $X$ is a nilpotent action.  By {\sc Sullivan} \cite[1971]{MR0494074} and {\sc Hilton-Mislin-Roitberg} \cite[1975]{MR0478146}, a nilpotent  space $X$ admits a {\em $P$-localization}  $\ell_X \colon X \to X_P$ which is a map inducing $P$-localization on   homotopy groups. 
When $P = \{ p \}$ we write $X_{p}$ for the $p$-localization and  when $P$ is empty we write $X_\Q$ for the rationalization of $X.$

Under reasonable hypotheses on $X$ and $Y$, function spaces behave well with respect to localization.   Hilton-Mislin-Roitberg  proved that if $X$ is a finite CW complex and $Y$ is a nilpotent space then the path-components of $\map(X, Y)$ are nilpotent spaces.  The components of $\map(X, Y)$ are of CW type in 
this case by {\sc Milnor} \cite[1959]{MR0100267}. Hilton-Mislin-Roitberg also proved that composition by $\ell_Y$ gives a $P$-localization map 
$$ (\ell_Y)_* \colon \map(X, Y; f) \to \map(X, Y_P; \ell_Y \circ f).$$
Below we  write $f_P = \ell_Y \circ f$.  {\sc Hilton-Mislin-Roitberg-Steiner} \cite[1978]{MR517093} 
obtained the same results if, alternately,  $X$ is a finite 
type CW complex and $Y$ is a nilpotent Postnikov piece. Here the CW type result is due to {\sc Kahn} \cite[1984]{MR733413}. {\sc M{\o}ller} \cite[1987]{MR902794} extended these results to the function space of relative liftings $$ \xymatrix{ A \ar[r]^u \ar[d]^i & E \ar[d]^p \\
B \ar@{.>}[ur] \ar[r]^f & X}$$
where here $i$ is closed cofibration and $p$ is a fibration.  In this case, $P$-localization is obtained by passing from $p \colon E \to X$ to its fibrewise $P$-localization
$p_{(P)} \colon E_{(P)} \to X$  as constructed by {\sc May} \cite[1980]{MR554323}. In particular,   fibrewise localization induces  $P$-localization $\Gamma(p;s) \to \Gamma(p_{(P)};s')$ for $X$ finite CW and  $F = p^{-1}(*)$ a nilpotent space.   {\sc M{\o}ller} \cite[1990]{MR1038731} proved the corresponding nilpotence and localization results for certain equivariant function spaces. 
{\sc Klein-Schochet-Smith} \cite[2009]{MR2574026}  extended the nilpotence and localization results from the case when $X$ is  finite CW to the case  $X$ is compact metric provided the corresponding function or section space is known {\em a priori} to be nilpotent.        

{\sc Bousfield-Kan} \cite[1972]{MR0365573} introduced a more general localization theory for subrings $R \subset \Q$   including homotopy completions $Y \to R_\infty Y. $   They proved that    $R$-completion (respectively,  $R$-localization)     induces $R$-completion (respectively,  $R$-localization) on the based function spaces  $\map_*(X, Y;0)$ when $X$ is a finite CW complex and $Y$ is nilpotent.   Further significant results on the behavior of function spaces under Bousfield-Kan localization and completion are discussed below.

\subsection{Rational Homotopy Theory of Function Spaces.}
  {\sc Quillen} \cite[1969]{MR0258031}  constructed an equivalence between   the homotopy category of simply connected rational CW complexes and   a homotopy category of connected, differential graded  Lie algebras (DGLAs) over $\Q$ initiating rational homotopy theory.  The Quillen minimal model of a simply connected space $X$ is a {\em minimal} DGLA $(\LL(X), d_X)$ which means 
  $\LL(X) = \L(V)$ is a free GLA and $d_X$ satisfies $d_X(V) \subseteq \left[\L(V), \L(V)\right].$  The rational homology and homotopy Lie algebra of $X$ are recovered via   isomorphisms
  $$ V \cong s^{-1}\widetilde{H}(X; \Q) \hbox{\, and \, } \pi_*(\Omega X) \otimes \Q, [\, ,\, ] \cong H_*(\LL(X)), [\, ,\, ].$$  
  
  {\sc Sullivan} \cite[1977]{MR0646078} constructed another categorical equivalence, here  between the homotopy theory of simply connected CW complexes and a homotopy category of connected differential graded algebras (DGAs) over $\Q$.   The Sullivan minimal model of a space $X$ is a minimal DGA $(\M(X), d_X)$ where $\M(X) = \Lambda V$
  is a free DGA with  $d_X(V) \subseteq \Lambda^+ V \cdot \Lambda^+ V$  with
  $$ V \cong \Hom(\pi_*(X), \Q) \hbox{\, and \, } H_*(\M(X)) \cong H^*(X; \Q).$$
More generally, a Sullivan model $(\A(X), d)$ for $X$ is a DGA admitting a chain equivalence to the de Rahm complex of  rational PL  forms on $X$.  In particular,
$H_*(\A(X)) \cong H^*(X; \Q)$ and $(\A(X), d) \simeq (\M(X), d_X)$ are homotopy equivalent in Sullivan's DGA category.  
  Comprehensive treatments of the subject were given by {\sc Tanr{\' e}} \cite[1983]{MR764769}  and 
{\sc F{\' e}lix-Halperin-Thomas} \cite[2001]{MR1802847}. 

Sullivan  described separate  models  for   the general path-component of a function space, $\map(X, Y;f)$, the space of self-equivalences $\aut(X)_\circ$ and the free loop space $\Lo X$ each within his framework of  DGAs.  We discuss these models and  their extensions and applications  now.

 \subsubsection{General Components.} Following the sketch by Sullivan, {\sc Haefliger} \cite[1982]{MR667163} constructed a (non-minimal) Sullivan model for the rational homotopy type of $\map(X, Y;f)$ where $f \colon X \to Y$ is a map of nilpotent spaces with $X$ finite.  The construction builds on the ideas of  Thom, described above. Let $p_r \colon Y_r \to Y_{r-1}$ with fibre $K(G_r, n_r)$ be a term in the principal  refinement of the Postnikov tower of $Y_\Q$ with $k$ invariant $k_{r-1} \colon Y_{r-1} \to K(G_r, n_r+1).$   
We then obtain  a pullback diagram $$\xymatrix{    \map(X, Y_r; (f_\Q)_r) \ar[rrr] \ar[d]_{(p_r)_*} &&& P\map(X, K(G_r, n_r +1);0) \ar[d] \\
   \map(X, Y_{r-1}; (f_\Q)_{r-1}) \ar[rrr]_{\! \! ( k_{r-1}\circ (f_\Q)_{r-1})_*} &&& \map(X, K(G_r, n_r +1);0) }$$ where the right fibration is  the   path/loop fibration.     Let $V = \bigoplus_r \Hom(G_r, \Q)$. 
 Since $X$ is finite, $X$ admits a finite   model $(A, d)$.   Write $\underline{A} = \Hom(A, \Q)$ for the dual to $A$ and grade $\underline{A}$ in negative degrees.    Thom's  calculation of $\pi_q(\map(X, K(G, n); 0))$ is then reflected in the grading on
the ordinary tensor product  $  \underline{A} \otimes V.$  
 Let $I$ denote the ideal of $\Lambda (\underline{A} \otimes V)$ generated by elements of degree $\leq 0.$ Haefliger described a differential $d$  on $\Lambda (\underline{A} \otimes V)/I$ in terms of the ``$k$-invariants''  $(k_{r-1}\circ (f_\Q)_{r-1})_*$  above   and proved 
 the result
    is a Sullivan model for $\map(X, Y;f).$ 

{\sc Bousfield-Peterson-Smith} \cite[1989]{MR989883} gave   an alternate construction, motivated by seminal work of {\sc Lannes} \cite[1987]{MR932261} in  $p$-local homotopy theory, discussed below.  The construction makes use of the    fact that $\map(X, \_\_)$ defines a functor on topological spaces that is right adjoint to the product functor $\_\_ \times X.$  In the category of DGAs, this   corresponds to the fact that $\Hom(A, \_\_)$ is right adjoint to the tensor product functor
$\_\_ \otimes A.$ In this  setting,  $\_\_ \otimes A$ has an   left adjoint, as well,      provided $A$ is finite. The construction is a version of   Lannes'  $T$-functor.   It is  conveniently written here as  $( \_\_  \, \colon A)$. Given a map $\psi \colon B \to A$ define $(B \, \colon A)_\psi$ to be the connected DGA determined by $\psi.$  Assume $X$ and $Y$ are nilpotent spaces with $X$ finite. 
Let $A$ be a finite Sullivan model for $X$. Bousfield-Peterson-Smith proved 
$ (\M(Y) \, \colon   A)_\psi $  is a Sullivan model for $\map(X, Y;f )$
where $\psi \colon \M(Y) \to A$ is a Sullivan model for $f \colon X \to Y.$

{\sc Brown-Szczarba} \cite[1998]{MR1407482} and   \cite[1998]{MR1473920} expanded on the work of Haefliger and Bousfield-Peterson-Smith. They constructed a model $(\Lambda W, d)$ for $\map(X, Y;f )$ where 
$$ W^q  = \left(\sum_n \pi_n(Y) \otimes H_{n-q}(X; \Q) \right)/ K^q$$ for certain  subspaces $K^q$.  The differential $d$ was described explicitly  in terms of the     coproduct on $H_*(X; \Q)$.  Here $X$ is a finite CW complex and $Y$ is nilpotent.  They deduced descriptions of the rational homotopy groups of $\map(X, Y;f).$ When $f$ is trivial they proved $K_q =0$ thus obtaining an  isomorphism of graded spaces 
$$\pi_*(\map(X, Y; 0)) \cong \left(H_*(X; \Q) \right)  \otimes \left( \pi_*(Y) \otimes \Q\right).$$
Again, the space $H_*(X; \Q)$ is assumed to be negatively graded. This last result was earlier proved by {\sc Smith} \cite[1994]{MR1225575}.

Applications of the Haefliger model include the following results: {\sc Vigu{\' e}-Poirrier} \cite[1986]{MR850369}  identified the rational homotopy Lie algebra of $ \map(X, Y; 0)$  for $X$  nilpotent of dimension  strictly less than the degree of the first nontrivial homotopy group of $Y$ via an isomorphism   $$ \pi_*(\Omega \map(X, Y; 0)) \otimes \Q, [ \, , \, ] \cong \left(H^*(X;\Q)\right) \otimes\left( \pi_*(\Omega Y)\otimes \Q)\right), \, [ \, , \, ]. $$ Here $H^*(X; \Q)$ is negatively graded and  the tensor product has the  GLA structure induced by the product and bracket on the terms. 
{\sc M{\o}ller-Raussen} \cite[1986]{MR808750} studied the rational homotopy classification problem for components 
of $\map(X, Y)$ with $Y = S^n, \C P^n$ for $X$ nilpotent and suitably rationally co-connected.  They obtained  complete classifications in these cases including descriptions of the rational homotopy types.  {\sc F{\' e}lix} \cite[1990]{MR1078101} proved the rational L.S. category of $\map(X, Y; 0)$  is often  infinite.  {\sc Smith} \cite[1997]{MR1471062}  studied the rational homotopy classification problem   for  $\map(G_1/T_1, G_2/T_2)$, where $G_1, G_2$ are classical compact Lie groups and $T_1, T_2$  maximal tori,  identifying   the rational type of certain components as  generalized flag manifolds.  {\sc Smith} \cite[1999]{MR1705124} gave an explicit description of  the Haefliger model     for $X$ and $Y$ {\em elliptic} spaces (simply connected spaces  having finite-dimensional rational homotopy and homology) with evenly graded rational cohomology   obtaining    examples of components of $\map(X, Y)$ of finite  L.S. category.
  
 {\sc Kotani}  \cite[2004]{MR2084591} used the Brown-Szczarba model to give necessary and sufficient conditions for the space $\map_*(X, Y)$ to be a rational $H$-space for $X$ a formal, nilpotent CW complex  of dimension $\leq$ the connectivity of $Y.$   {\sc F{\' e}lix-Tanr{\' e}} \cite[2005]{MR2153111} generalized this result replacing the formality of $X$ by a condition involving L.S. category.   {\sc Buijs-Murillo} \cite[2006]{MR2442961}  constructed the Brown-Szczarba model within the simplicial category framework for   rational homotopy theory due to  {\sc Bousfield-Gugenheim} \cite[1976]{MR0425956}.   They obtained a functorial version of the Brown-Szczarba model in this setting  and used this model, in {\sc Buijs-Murillo}  \cite[2008]{MR2442961},   to  identify  the rational homotopy Lie algebra of components $\map(X, Y;f)$ with $X,Y$ restricted, as usual, to nilpotent spaces with $X$ finite.    {\sc Kuribayashi-Yamaguchi} \cite[2006]{MR2220679} combined the Haefliger and Brown-Sczarba approaches  to obtain a rational splitting   of   $\map_*(X \cup_\alpha e^{k+1},Y; 0)$ where $\alpha$ is an   attaching map. Under certain restrictions on $X$, $Y$ and $\alpha$ they proved
$$\map_*(X \cup_\alpha e^{k+1},Y; 0) \simeq_\Q \map_*(X, Y;0) \times \Omega^{k+1}Y.$$
{\sc Hirato-Kuribayashi-Oda} \cite[2008]{MR2369043} applied the Brown-Szczarba model to the study of the rational evaluation subgroups, i.e.,  the image of the map induced on  rational homotopy groups by the evaluation map $\omega_f \colon \map(X, Y; f) \to Y.$ 
{\sc Buijs-F{\' e}lix-Murillo} \cite[2009]{MR2531370} used the Brown-Szczarba model to study the rational homotopy type of the homotopy fixed point of a circle action.  

The higher rational homotopy groups of $ \map(X, Y;f )$   for suitable $X$ and $Y$ can be described directly in terms the homology of certain    DG space of  derivations.   In the DGA setting,  given a map $\psi \colon (A,d) \to (B, d)$ of DGAs  let $\Der_n(A, B; \psi)$ denote the space of linear maps $\theta \colon A^* \to B^{*-n}$ satisfying, for $x, y \in A,$ 
the identity $\theta(xy) = \theta(x) \psi(y)+ (-1)^{n|x|}\psi(x)\theta(y).$ 
In the DGLA setting, let $\Der_*(L, K; \psi)$ denote the space of degree raising linear maps satisfying the corresponding identity. In both cases, a degree $-1$ differential is given by
  $D(\theta) = d  \circ \theta - (-1)^{n} \theta \circ d .$    When $X$ and $Y$ nilpotent spaces with $X$ finite   $$\pi_n(\map(X, Y;f)) \cong H_n(\Der(\M(Y), \M(X); \M(f))).$$
for $n \geq 2$.   This result is due to Sullivan for the case $f = 1$ as discussed below. The general result was proved, independently, by {\sc Block-Lazarev} \cite [2005]{MR2132759}, {\sc Lupton-Smith} \cite[2007]{MR2292124} and {\sc Buijs-Murillo}
\cite[2008]{MR2442961}.   The rationalization of the  fundamental group 
$\pi_1(\map(X, Y;f)) $ is, in general, nonabelian. {\sc Lupton-Smith} \cite[2007]{MR2302588}   proved the rank of $\pi_1(\map(X, Y;f))_\Q$ is  the dimension of $H_1(\Der(\M(Y), \M(X); \M(f))).$ {\sc Buijs-Murillo}
\cite[2008]{MR2442961} extended this to an identification of the  Malc'ev completion of
$\pi_1(\map(X, Y;f))_\Q$.

Within the framework of Quillen minimal models,   we have   an  isomorphism   
$$\pi_n(\map(X, Y; f))\otimes \Q \cong H_n(\Rel(\ad_{\LL(f)}))$$
for $n \geq 2$ for $X$ and $Y$ simply connected CW with $X$ finite.  Here    $\Rel_*(\ad_{\LL(f)})$
is the mapping cone of the chain map $$\ad_{\LL(f)} \colon \LL(Y) \to  \Der_*(\LL(X), \LL(Y); \LL(f))$$ given by $\ad_{\LL(f)}(y)(x)  = [  \mathcal{L}(f)(x), y]$ for $x \in \LL(X), y \in \LL(Y).$
This result was proved for the identity component by {\sc Tanr{\' e}} \cite[1983]{MR764769} and {\sc Schlessinger-Stasheff} \cite[preprint]{S-S}. The result for the general component was proved by {\sc Lupton-Smith} \cite[2007]{MR2292125}.
{\sc Lupton-Smith} \cite[2010]{preprintLS} identified rational Whitehead products in terms of this  identification.   
{\sc Buijs-F{\' e}lix-Murillo} \cite[2009]{MR2515825} described  a Quillen model
for function spaces  and obtained a result on the exponential growth of rational homotopy
groups of function spaces. 

The   homotopy classification problem    for gauge groups corresponding to principal $G$-bundles $P \colon E \to X$  is trivial after rationalization for $X$ finite CW and $G$ a compact Lie group. In this case, $BG$ is a rational $H$-space.  As mentioned above, {\sc Crabb-Sutherland}
\cite[2000]{MR1781154} used this fact to prove the fibrewise localization of the  universal $G$-adjoint bundle $\Ad(P_G) \colon E_G \times_G G^{\ad}$  is equivariantly trivial. Their result implies a rational $H$-equivalence $$(\G(P)_\circ)_\Q \simeq_H \map(X, G_\Q; 0).$$
The rational homotopy groups of $\G(P)$ may thus be computed as $$\pi_*(\G(P)) \otimes \Q  \cong H_*(X; \Q) \otimes (\pi_*(G) \otimes \Q)$$
since $G_\Q$ is   a product of Eilenberg-Mac \! Lane spaces where here, as above,   $H_*(X; \Q)$ is negatively graded. 

{\sc Lupton-Phillips-Schochet-Smith} \cite[2009]{MR2439407} proved a related result in the context of commutative Banach algebras. Let $A$ be unital, commutative Banach algebra and 
$\mathrm{GL}_n(A)$ the group of $n \times n$ invertible matrices with coefficients in $A.$ Then
$$ \mathrm{GL}_n(A)_\circ \simeq_\Q \prod K(V_n , n) \hbox{\,  where \, } V = \Check{H}_*(\mathrm{Max}(A); \Q) \otimes \Lambda(s_1, \ldots, s_{2n-1}).$$
Here $\mathrm{Max}(A)$ is the maximal ideals space and $s_{2i-1}$ is of degree $2i-1$. 
{\sc Klein-Schochet-Smith} \cite[2009]{MR2574026} extended this  result to  the  group of unitaries $UA_\zeta$ where $A_\zeta$ is the $C^*$-algebra sections of a complex $n$-matrix bundle $\zeta$ over a compact metric space $X.$  The latter result is based on an extension of the result $(\G(P)_\circ)_\Q \simeq_H  \map(X, G_\Q; 0)$  from the case $X$ finite CW to the case $X$ compact metric. 

\subsubsection{Spaces of Self-Equivalences.} The rational homotopy type of a connected grouplike space $G$  is completely determined by isomorphism type of the Samelson Lie algebra $\pi_*(G) \otimes \Q, [ \, , \, ]$ (c.f. {\sc Scheerer} \cite[1985]{MR788673}).  {\sc Sullivan} \cite[1977]{MR0646078} identified the rational Samelson Lie algebra of the space $\aut(X)_\circ$ for $X$ a simply connected finite CW complex via an  isomorphism:  
 $$ \pi_*(\aut(X)_\circ) \otimes \Q,   [ \, , \, ] \cong H_*(\Der(\M(X))), [ \, , \, ].$$
  A corresponding identity in Quillen's  DGLA framework for rational homotopy theory was given by {\sc Tanr{\' e}} \cite[1983]{MR764769} and {\sc Schlessinger-Stasheff} \cite[preprint]{S-S}: 
$$ \pi_*(\aut(X)_\circ) \otimes \Q, [ \, , \, ] \cong H_*(\Rel(\ad_{\LL(X)})), [ \, , \, ].$$
Here the bracket on the mapping cone of $\ad_{\LL(X)}$ is induced from that on 
$\mathcal{L}(X)$. 

 Sullivan's identity  connects the monoid $\aut(X)$ to a fundamental open conjecture in rational homotopy theory.  Let $X$ be a simply connected elliptic
CW complex with evenly graded rational cohomology. We refer to such spaces as {\em $F_0$-spaces.} The class includes (products of) spheres, complex projective spaces and homogeneous spaces $G/H$ with $G$ a compact Lie group and $H \subset G$ a closed subgroup of maximal rank.  Motivated by this last case, {\sc Halperin} \cite[1978]{MR0515558} conjectured that the rational Serre spectral sequence collapses at the $E_2$ term for all orientable fibrations with fibre  an $F_0$-space. {\sc Thomas} \cite[1981]{MR638617} and {\sc Meier} \cite[1981]{MR649203} independently proved that Halperin's conjecture   is equivalent
to the condition $H_{\mathrm{even}}(\Der(\M(X))) = 0$ for an $F_0$-space $X$. Thus, by Sullivan's identity, Halperin's conjecture holds for $X$ if and only if  $\aut(X)_\circ$ is rationally  equivalent to a product of odd spheres.

 The Halperin conjecture has been affirmed in several special cases   including for K{\" a}hler manifolds by Meier, for homogeneous spaces of maximal rank pairs by {\sc Shiga-Tezuka} \cite[1987]{MR894562} 
 and     for $F_0$-spaces with rational cohomology generated by $\leq 3$ generators by {\sc Lupton} \cite[1990]{MR1082989}.  {\sc Hauschild}  \cite[1993]{MR1226781} \cite[2001]{MR1817009} computed rational homotopy groups of $\aut(X)$     for $X$ an $F_0$-space and  of $\aut(p)$ for $p$ a fibration with fibre $X$.    {\sc F{\' e}lix-Thomas}  \cite[1994]{MR1297839} studied  the rational homotopy of $\aut(X)$ for various homogeneous spaces giving   explicit calculations. {\sc Grivel} \cite[1994]{MR1255925} proved that for, $X$ an $F_0$-space, 
 $$H_{\mathrm{even}}(\Der(\M(X))) \cong \Der_{\mathrm{even}}(H(\M(X)))$$
 and used this result to give a formula for $\pi_{\mathrm{odd}}(\aut(X)_\circ) \otimes \Q.$ 
 
 {\sc Salvatore} \cite[1997]{MR1443426} proved  the nilpotency of the Lie algebra $H_*(\Der( \M(X)))$ coincides   with the {\em rational  homotopical nilpotency} of the monoid $\aut(X)_\circ $ ---  the least integer $n$ such that the $n$-fold commutator for $\aut(X)_\circ$ is rationally trivial.  He    calculated the rational homotopical nilpotency of $\aut(X)_\circ$ for   $X$
    a rational two-stage Postnikov system and proved the monoid $\aut(S^{2n-1} \vee S^{2n-1})_\circ$ is not rationally homotopy nilpotent.  {\sc Smith} \cite[2001]{MR1836007} computed the rational homotopy nilpotency of $\aut(X)_\circ$ for certain spaces  $X$ admitting a two-stage Sullivan model.  Recently, {\sc F{\' e}lix-Lupton-Smith} \cite[preprint]{preprintFLS} obtained a formula, in the spirit of Sullivan's above, for the rational Samelson Lie algebra of the monoid $\aut(p)_\circ$ of fibre-homotopy self-equivalences of a fibration $p \colon E \to B$ of simply connected finite CW complexes:
 $$\pi_*(\aut(p)_\circ) \otimes \Q, [ \, , \, ] \cong H_*(\Der_{\Lambda V}(\Lambda W \otimes \Lambda V)), [ \, , \, ].$$
Here $(\Lambda W, d_B) \to (\Lambda W \otimes  \Lambda V, D)$ is the Koszul-Sullivan model of the fibration and $\Der_{\Lambda V}(\Lambda W \otimes \Lambda V)$ denotes the DGLA of derivations vanishing on $\Lambda V.$ 
 
As for the rational homotopy of Gottlieb group $G_n(X)$ and the evaluation map $\omega \colon \aut(X)_\circ \to X,$  {\sc Lang} \cite[1975]{MR0367986} proved $G_*(X)_\Q \cong G_*(X_\Q)$ for $X$ a finite simply connected CW complex.   {\sc F{\' e}lix-Halperin} \cite[1982]{MR664027} identified the
 rationalized Gottlieb groups $G_n(X) \otimes \Q$  for these $X$   in terms of the Sullivan identification
$\pi_*(X) \otimes \Q  \cong V$ where $\M(X) = \Lambda V.$  Here an element $v \in V_n$ corresponds to a rational Gottlieb element if the dual map $v \mapsto 1$ extends to a  derivation cycle in $\Der_n(\M(X)).$ They used this result to prove two global results on the rationalized Gottlieb groups of a simply connected CW complex $X$ of finite rational L. S. category: $$G_{\mathrm{even}}(X) \otimes \Q = 0
\hbox{\,   and \, dim}(G_{\mathrm{odd}}(X) \otimes \Q) \leq \mathrm{cat}_\Q(X).$$ 
A Quillen model description of the rationalized Gottlieb group was given by {\sc Tanr{\' e}} \cite[1983]{MR764769}:
$$G_n(X)_\Q \cong \ker\{ H(\ad_{\LL(X)}) \colon \LL(X) \to \Der(\LL(X)) \}.$$
 Rationalized Gottlieb groups have been calculated by many authors using various means including {\sc Smith} \cite[1996]{MR1381587}, {\sc Lupton-Smith} \cite[2007]{MR2292124} and \cite[2007]{MR2292125}, {\sc Hirato-Kuribayashi-Oda} \cite[2008]{MR2369043} and {\sc Yamaguchi} \cite[2008]{MR2475490}. 
 {\sc F{\' e}lix-Lupton} \cite[2007]{MR2337558} proved the evaluation map $\omega \colon \aut(X)_\circ \to X$  is rationally homotopy trivial if and only if it is trivial  on rational homotopy groups for $X$ a finite, simply connected CW complex. 

\subsubsection{The Free Loop Space.} {\sc Vigu{\'e}-Poirrier-Sullivan} \cite[1976]{MR0455028} 
constructed a Sullivan model for $\Lo X$ when $X$ is a simply connected CW complex.  Let $(\Lambda V, d)$ be the minimal model for $X$. Their model for $\Lo X$ is given by $$(\Lambda V \otimes \Lambda  sV ,\delta) \hbox{\, with \, } \delta(v) = dv \hbox{ \, and \, } \delta(sv) = -sd(v) \hbox{\, for $v \in V$}$$
where $s$ is the degree $-1$ derivation  of $\Lambda V \otimes \Lambda sV$  defined by setting $s(v) = sv$
and $s(sV) = 0.$   {\sc Halperin} \cite[1981]{MR633568}  constructed  a related model in the non-simply-connected case under restrictions on the component.    
As discussed above,  Vigu{\' e}-Poirrier-Sullivan used  their model to prove the Betti numbers of $LX$ are unbounded if $H^*(X; \Q)$ requires at least two generators.  In fact, 
 they showed the Betti numbers of $\Lo  X$ grow exponentially in this case.
{\sc Vigu{\' e}-Poirrier} \cite[1984]{MR777376} proved that the same  is true for wedges of spheres and manifolds of L.S. category less than $2$. She conjectured the Betti numbers of $\Lo X$ grow exponentially   for all finite, simply connected $X$ with infinite dimensional rational homotopy.  She   gave examples of spaces $X$ with finite-dimensional rational homotopy for which 
$$\limsup\frac{\log\sum^n_{i=1}\beta_i}{\log n}=\dim(\pi_{{\rm odd}} (X)\otimes \Q).$$
  {\sc Lambrechts} \cite[2001]{MR1834084} affirmed Vigu{\' e}-Poirrier's conjecture for  the class of simply connected, coformal, finite complexes.   

{\sc Dupont-Vigu{\' e}-Poirrier} \cite[1998]{MR1651384} proved a basic result concerning the question of formality for the free loop space.  Given  a  simply connected  CW complex $X$ with Noetherian rational cohomology,  the space $\Lo X$ is formal if and only if $X$ is a rational $H$-space.  {\sc Yamaguchi} \cite[2000]{MR1772441}  generalized this to the function space $\map(X, Y;0)$ for $X$ simply connected CW of dimension less than the connectivity of $Y$ an elliptic space. He showed that, again,   $\map(X, Y;0)$ is formal if and only if $Y$ is a rational $H$-space.  {\sc Vigu{\' e}-Poirrier} \cite[2007]{MR2326935}  proved a further result in this vein showing, with the same dimension/connectivity hypotheses,  that $\map(X, Y; 0)$ is formal  if and only if $Y$ is a rational $H$-space provided the odd rational Hurewicz homomorphism of $X$ is nontrivial. 

{\sc F{\' e}lix-Thomas-Vigu{\' e}-Poirrier} \cite[2007]{MR2283106} studied the string topology operations on $\H_*(\Lo M; \Q)$ within the framework of Sullivan minimal models. They gave a description of the loop-product
and the string-bracket in this setting making explicit computations.  They also proved the 
isomorphism of graded spaces due to {\sc Cohen-Jones} \cite[2002]{MR1942249}   
$$ \H_*(\Lo M; \Q) \cong HH^*(C^*(M), C^*(M); \Q)$$
is an isomorphism of Gerstenhaber algebras.  {\sc F{\' e}lix-Thomas} \cite[2008]{MR2415345}  extended this last result proving the above is an isomorphism of BV-algebras.

\subsection{ Function Spaces and $p$-Localization.} Function spaces are central    to the  theory of homotopy  localizations and  completions   with respect to subrings $R \subset \Q$.  As mentioned above, {\sc Bousfield-Kan} \cite[1972]{MR0365573} proved that their $R$-completion functor induces $R$-completion on the based function space $\map_*(X, Y; 0)$  for  $X$ finite CW and   $Y$ nilpotent.  This result  was used in the proof of the fracture lemma for $R$-completions.  A first reduction in Miller's proof of the Sullivan conjecture is a weak equivalence $$\map_*(X,Y; 0) \simeq_w \map_*(X, R_\infty Y; 0)$$ where $R = \mathbb{F}_p$ and $R_\infty$ is the $R$-completion functor.  Here $Y$ is a nilpotent space and $X$ is a connected, $\Z[\frac{1}{p}]$-acyclic space. 
 Function spaces also feature in the theory of     homotopy localization and cellularization. Given a map $f\colon X\to Y,$   a space $Z$ is defined to be $f$-local if the induced map of function spaces $f^* \colon \map_*(Y,Z)\to \map_*(X,Z)$ is a weak homotopy equivalence. {\sc Dror-Farjoun}    \cite[1996]{MR1392221} constructed  $f$-localizations  and showed the known localization functors all occur as  special cases, for suitable choices of the map $f$.    As usual, we consider work here which  focuses explicitly on the homotopy type of function spaces. 
 
\subsubsection{Maps out of a Classifying Space.}   The   function spaces  $\map(B\pi, X)$ and $\map_*(B\pi, X)$ for $\pi$ a finite group appear in major developments 
in homotopy theory in the $p$-local category. In celebrated work, {\sc Miller} \cite[1984]{MR750716}
affirmed the Sullivan conjecture, proving  $$\pi_n(\map_*(B\pi, X)) = 0 \hbox{\, for  all \, } n \geq 0,$$
for $\pi$ a locally finite group and $X$ a connected, finite CW complex.   Using the $R$-completion theorem for $\map_*(X, Y)$  mentioned above,  the problem
reduces to proving the weak triviality of spaces $\map_*(B\Z_p, R_\infty Y)$ where $R_\infty$ is the Bousfield-Kan $p$-completion functor. Miller proved  the latter fact by establishing the vanishing of certain  Ext-sets  in a  category  
of unstable  modules over the mod $p$ Steenrod algebra.    

Miller's Theorem had many important consequences for function spaces.  {\sc Zabrodsky} \cite[1991]{MR1177330} connected the result to the study of phantom maps. He also    obtained  the following extension. Let $W$ be  a connected CW complex with finitely many, locally finite homotopy groups. Then 
$$\pi_n(\map_*(W, X)) = 0 \hbox{\, for  all \, } n \geq 0,$$
for $X$ any connected, finite CW complex.  Zabrodsky also proved that, if $P \colon E \to B$ is a principal $G$-bundle with $\map_*(G, Y)$ contractible,  then $\map_*(B, Y) \to \map_*(E, Y)$ is a homotopy equivalence.  
Miller's result  had equivariant generalizations in terms of fixed point and   homotopy fixed point set.   In particular,  if $X$ is a $\pi$-space for $\pi$ a $p$-group then 
$$R_\infty(X^\pi) = (R_\infty X)^{h\pi}$$ where $R=\Z/p$ and $X^\pi$ is the fixed point set while $X^{h\pi}$ is the homotopy fixed point sets (c.f. {\sc Carlsson} \cite[1991]{MR1091616}).    {\sc Dwyer-Zabrodsky} \cite[1986]{MR928826}  applied this latter result to obtain a mod $p$ decomposition
$$ \map(B\pi, BG) \simeq_{p}  \coprod_{\rho}BC(\rho).$$
Here $\pi$ is a finite $p$-group, $G$ is a compact Lie group and the disjoint union is over
$G$-conjugacy classes of homomorphisms $\rho \colon \pi \to G$. As usual,  $C(\rho)$ denotes the centralizer of the image of $\rho$ in $G$.
{\sc Friedlander-Mislin} \cite[1986]{MR827361} gave conditions on a Lie group $G$ such that $\map_*(BG, R_\infty X)$ is weakly trivial.  Here $X$ is nilpotent and $R_\infty$ is $p$-completion.  {\sc McGibbon} \cite[1996]{MR1328362} proved $\map_*(W, R_\infty X)$ is weakly contractible for $W$ a connected infinite loop space with torsion fundamental group. {\sc Strom} \cite[2005]{MR2029919} proved that 
 if $\map_*(X,S^n)$ is weakly contractible  for all sufficiently large $n$ then $\map_*(X, Y)$ is actually weakly contractible for any nilpotent, finite CW complex $Y$. He  thus obtained a  method for recognizing spaces $X$ satisfying the conclusion of Miller's Theorem.

 {\sc Lannes} \cite[1987]{MR932261} and \cite[1992]{MR1179079} complemented and extended Miller's work. Let $\mathcal{U}$ and $\mathcal{K}$ denote, respectively,  the category of unstable modules and  algebras  over   $\mathcal{A},$ the   mod $p$ Steenrod algebra. Lannes constructed the $T$-functor which is left-adjoint to the tensor product functor $ \_ \_ \otimes_{\mathcal{U}} H^*(B\Z_p; \Z_p)$ on $\mathcal{U}.$  He showed, among other properties, that    $T$ is exact,  preserves tensor products and restricts to a functor on $\mathcal{K}$.  These results are used to prove the natural map $$\Theta_X \colon  T(H^*(X; \Z_p)) \to H^*(\map(B\Z_p, X); \Z_p)$$
is an  isomorphism of unstable $\mathcal{A}$-algebras whenever $T(H^*(X; \Z_p))$ is trivial in degree $1$ and of finite type. 
 
 {\sc Aguad{\' e}} \cite[1989]{MR1005014} computed $T(M)$ for certain   $\mathcal{A}$-algebras  $M$ including subalgebras of $\Z_p[x_1, \ldots, x_n]$ invariant under the general linear action. Here $p$ is an odd prime. {\sc Dror-Smith} \cite[1990]{MR1098968} constructed an Eilenberg-Moore spectral sequence for computing $\Theta_X$ above and gave a geometric interpretation of the $T$-functor. {\sc Dwyer-Wilkerson}  \cite[1990]{MR1098969} proved if $\pi$ is a locally finite group and $X$ is a  simply connected $p$-complete space with $H^*(X; \Z_p)$ finitely generated as an algebra, then   $\map_*(B\pi,X; 0)$  is weakly contractible.  {\sc Kuhn-Winstead} \cite[1996]{MR1404914}  proved   $$\widetilde{H}^*(X; \Z_p) = 0 \, \,\Longrightarrow \, \, \widetilde{H}^*(\map(B\Z_p, X); \Z_p) = 0.$$
 More generally, they showed  the image of $\Theta_X$    consists of $\Z^{\wedge}_p$-integral classes if   $H^*(X; \Z_p)$ does where a $\Z^{\wedge}_p$-integral class in $H^*(\_ \_ ; \Z_p)$ is a class  in the reduction from $H^*(\_ \_ ;\Z^{\wedge}_p ).$
 {\sc Dehon-Lannes} \cite[2000]{MR1793415} proved that if $X$ is $p$-complete and $H^*(X; \Z_p)$ is Noetherian and generated in even degrees then    $H^*(\map(B\Z_p, X); \Z
 _p)$ is Noetherian and $\map(B\Z_p, X)$ is $p$-complete. {\sc Aguad{\' e}-Broto-Saumell} 
\cite[2004]{MR2039756} introduced the notion of $T$-representability for a space $X$  and proved it a sufficient condition for $\Theta_X$ to be an isomorphism. They gave an example of a $p$-complete space $X$ for which $\Theta_X$ is not an isomorphism.  

We mention   some further results falling under the current heading. 
{\sc Dwyer-Mislin} \cite[1987]{MR928824} identified the homotopy type of the nontrivial components of $\map_*(BS^3, BS^3)$. The null component is contractible by {\sc Zabrodsky}  \cite[1991]{MR1177330}, as  mentioned above.  {\sc Jackowski-McClure-Oliver} \cite[1992]{MR1147962} determined the homology with coefficients in a finite group  of the components of    $\map(BG, BG)$ for $G$ a simple compact, connected Lie group. 
Building on these results, {\sc Andersen-Grodal}  \cite[2009]{MR2476779} expressed the classification of $p$-compact groups  with component group a $p$-group $\pi$  in terms of the homotopy types of the components of a space of maps out of  $B\pi.$

{\sc Blanc-Notbohm} \cite[1993]{MR1112487} proved 
$$\map(BG, BH; f )_p^\wedge \simeq \map(BG, (BH)_p^\wedge; f_p^\wedge) $$
for $G, H$ compact Lie groups. 
{\sc Broto-Levi} \cite[2002]{MR1759888} proved     $$\aut((B\pi)^{\wedge}_p)_\circ \simeq K(Z(\pi/O_{p'}(\pi)), 1)$$ for $\pi$   a finite group. Here  $O_{p'}(\pi)$ denotes the  maximal normal $p'$-subgroup of $\pi$. The proof uses the Bousfield-Kan spectral sequence to prove asphericality and  the $Z^*$-theorem from group theory to compute the fundamental group. 
\subsubsection{Algebraic Models.} In this final section, we discuss results on modeling the $p$-local homotopy theory of function spaces. 

 {\sc Dwyer} \cite[1979]{MR551014}    proved  that  a version of Quillen's rational homotopy theory extends to the  homotopy theory of {\em tame spaces} $X$ which are $(r-1)$-connected CW complexes, $r\geq 3$, with  $\pi_{r+k}(X)$  uniquely $p$-divisible for  all primes $p$ with $2p-3\leq k$.
He  proved the homotopy category of tame spaces endowed with an appropriate model structure is   equivalent to a homotopy category of $(r-1)$-reduced integral DGLAs. 
Tame homotopy theory admits a   version   in the spirit of Sullivan's approach to rational homotopy theory by {\sc Cenkl-Porter} \cite[1983]{MR700985}.   {\sc Scheerer-Tanr{\' e}} \cite[1988]{MR972081}  identified homotopy invariants, e.g., homology and the homotopy Lie algebra, in Dwyer's framework.  {\sc Anick} \cite[1989]{MR991015} and 
\cite[1990]{MR1048177} gave an alternate approach    using a classical construction of   {\sc Adams-Hilton} \cite[1955]{MR0077929} on the Pontryagin algebra $H_*(\Omega X; R)$ for $R \subset \Q$.  He constructed a  DGLA  model  over $R$  for  the category   $\mathrm{CW}^m_r$  consisting of $(r-1)$-connected complexes of dimension $m$ when all primes $p$ with $m >  pr$ are invertible in $R$.   

The description of function spaces in tame homotopy theory has been undertaken in several works. {\sc Anick-Dror-Farjoun} \cite[1990]{MR1098963} described a  simplicial  skeleton  of the space $\map_*(X, Y)$ for $R$-local spaces $X, Y  \in \mathrm{CW}^m_r$. Given suitably reduced   DGLAs $L$ and $K$ over $R$ they constructed a function complex $\Hom(L, K)$   giving an explicit description  through a range of dimensions.
{\sc Scheerer-Tanr{\' e}} \cite[1992]{MR1170636}  described $R$-local  homotopy theory in a suitable category of DG coalgebras of $R$. 
They constructed an adjoint to the wedge functor  in this context and used it to give a model for the space $\map_*(X, Y)$ for $R$-local $X, Y \in \mathrm{CW}^m_r$.  As an example, they described a  model for the $R$-localization of $\map_*({\mathbb H}P^2,M_t)$, where $M_t$ denotes a tamed Moore space and $R$ is suitably chosen. {\sc F{\' e}lix-Thomas} \cite[1993]{MR1204401} obtained a $p$-local  decomposition 
$$ \map_*(\Sigma X, Y) \simeq_p \prod_{i=2}^{n}(\Omega^{i+1}Y)^{\beta_i(X)}$$
for $X$ simply connected with torsion-free homology of dimension $n< 2p$ and $Y$ $(r-1)$-connected with $r > n+1.$   {\sc Scheerer} \cite[1994]{MR1255936} recovered this result as a special case of a corresponding decomposition for  $\map_*(C, Y)$ for $C$ a co-H-space.  
{\sc Dupont-Hess} \cite[1999]{MR1704545},  \cite[2002]{MR1918186} and  \cite[2003]{MR1981874} used Anick's framework   to obtain a model for the mod  $p$  cohomology   of the free loop space $\Lo X$ for a simply connected space $X \in \mathrm{CW}_r^m$ and  prime $p$ with $m \leq pr$.  They constructed a DGA over $\Z_p$, built from Anick's model,  and proved the homology  of this algebra is isomorphic   to $H^{*}(\Lo X;\Z_p)$. 

{\sc Dwyer-Kan} \cite[1980]{MR584566}  identified   function complexes in a simplicial homotopy theory    category $L^H(\bf{M})$ of a general Quillen model  category ${\bf M}$. Here $L^H$ is their  ``hammock'' localization functor.  They showed   $L^H({\bf M}(X, Y))$ for $X, Y \in {\bf M}$ 
behaves properly with respect to simplicial and cosimplicial resolutions. This approach yielded, in particular, a good model for the monoid of self-equivalences for arbitrary model categories.  In {\sc Dwyer-Kan} \cite[1983]{MR705421}, they identified the homotopy type of the function complex as a homotopy inverse limit for $X$ cofibrant and $Y$ fibrant.

Recently, {\sc Fresse} \cite[preprint]{Fressepreprint} 
 constructed an algebraic model for the  homotopy type of $\map(X, Y)$  for $X$ a finite complex and $\pi_n(Y)$ a finite $p$-group for all $n.$  Let $N^*(Y; \overline{\Z}_p)$ denote the normalized cochain complex where $\overline{\Z}_p$ is the algebraic closure of the field of $p$-elements.  Then $N^*(Y;\overline{\Z}_p)$ is an $E_\infty$-algebra and, more generally, an algebra over the Barrett-Eccles operad $\mathcal{E}.$ 
{\sc Mandell} \cite[2001]{MR1791268} proved this algebra structure determines the homotopy type of $Y.$
Fresse constructed a Lannes $T$-functor left adjoint to the tensor product in the category of algebras over $\mathcal{E}.$ As a consequence, he obtained  a model
$(N^*(Y; \overline{\Z}_p) \,  \colon N^*(X; \overline{\Z}_p))$ for $\map(X, Y).$  {\sc Chataur-Kuribayashi} \cite[2007]{MR2320001}  extended this result to the case $Y$ is connected and nilpotent of finite type and made  calculations with the resulting spectral sequence. 
 
  \nocite{*}
\providecommand{\MR}{\relax\ifhmode\unskip\space\fi MR }
\bibliographystyle{amsplain}
\def\cprime{$'$}
\providecommand{\bysame}{\leavevmode\hbox to3em{\hrulefill}\thinspace}
\providecommand{\MR}{\relax\ifhmode\unskip\space\fi MR }
\providecommand{\MRhref}[2]{%
  \href{http://www.ams.org/mathscinet-getitem?mr=#1}{#2}
}
\providecommand{\href}[2]{#2}

\end{document}